\documentclass{amsart}
\usepackage{graphicx} 
\usepackage{biblatex}
\usepackage{style}
\title{Towards 2-derivators for formal $\infty$-category theory}
\author{Nicola Di Vittorio}
\date{}

\begin{document}
\begin{abstract}


Derivators, introduced independently by Grothendieck and Heller in the 1980s, provide a categorical framework for studying homotopy theory. They are based on the idea that, while the homotopy 1-category of a single model category or $(\infty, 1)$-category retains only limited information, the structured collection of homotopy 1-categories of diagram categories often suffices for many homotopical purposes. In this paper, we introduce a set of axioms for a 2-dimensional analog of derivators: a refinement of the homotopy 2-category of an enriched model category or $(\infty, 2)$-category into a coherent system of homotopy 2-categories of higher categories of diagrams. 
We show that these axioms are satisfied in a variety of models, including standard ones related to
$(\infty, 1)$-category theory. Moreover, we prove that the axioms are preserved under a certain shift operation.
\end{abstract}

\maketitle
\epigraph{Caelum, non animum mutant qui trans mare currunt.}{\textit{Horace}}
\section{Introduction}\label{chap:intro}

In recent years, there has been a paradigm shift in several areas of Mathematics from the use of algebraic structures to that of homotopy coherent structures, whose axioms do not hold strictly but rather are part of extra data satisfying higher dimensional coherences \cite{riehl2018homotopy}. Analogously to how algebraic structures can be studied using categories, homotopy coherent structures are better understood via $(\infty,1)$-categories (or $\infty$-categories for short). To put it simply, $\infty$-categories are category-like structures having not only objects and morphisms but also $k$-dimensional morphisms, or \emph{k-cells}, for every natural number $k$. These $k$-cells are required to be invertible in a suitable sense for every $k>1$. This intuition has been made precise using a variety of different models like simplicial categories \cite{bergner2007model}, quasi-categories \cite{boardman2006homotopy,joyal2008theory, luriehigher}, complete Segal spaces \cite{rezk2001model}, and so on. To name a few places where $\infty$-categories are commonly found, Derived Algebraic Geometry \cite{lurie2004derived} uses $\infty$-topoi $-$ a special kind of $\infty$-category that is a generalisation of Grothendieck topoi from ordinary category theory $-$ as a framework in which to study sheaves of homotopy types. The stable $\infty$-category of spectra \cite{luriehigher} provides a very well-behaved environment for Stable Homotopy Theory. In \cite{blumberg2013universal}, the authors characterise Higher Algebraic K-theory by means of a universal property that is expressed in terms of $\infty$-categories and the list goes on.

Just as (small) categories are objects of the 2-category $\Cat$, (small) $\infty$-categories assemble into an $(\infty,2)$-category $\Cat_{\infty}$. Roughly speaking, $(\infty,2)$-categories resemble $\infty$-categories in the sense that they have $k$-cells for every natural $k$, with the difference that the $k$-cells are invertible for $k>2$. Fixing a natural number $n$, one can work with higher categories in which all $k$-dimensional morphisms are invertible for $k>n$; these are called $(\infty,n)$-categories. They play an important role in the study of Topological Quantum Field Theories as witnessed by Lurie's work \cite{lurie2008classification} on the Cobordism Hypothesis \cite{baez1995higher}. The theory of $(\infty,2)$-categories has also been used as a setting for Geometric Langlands in Gaitsgory and Rozenblyum's books \cite{gaitsgory2019study,gaitsgory2017study}. Applications aside, $(\infty,2)$-categories are a natural setting for the theory of $\infty$-categories. In the spirit of papers such as \cite{LackS:fortm, Street:fib-yon, StreetR:fortm}, that develop aspects of the theory of categories within a 2-category with enough structure, one can seek to do the same for $\infty$-categories in a nice enough $(\infty,2)$-category. A major advantage of this approach is that the theory of $\infty$-categories can be developed without having to actually define $\infty$-categories explicitly, treating them as a primitive notion instead. Unlike categories, whose definition is reasonably straightforward, making the concept of $\infty$-categories precise takes a considerable amount of work and different models are useful in different occasions. However, thinking about $\infty$-categories as objects inside an $(\infty,2)$-category comes at the cost of defining what an $(\infty,2)$-category is to begin with. As one might expect, there are many models for them. Even though all of the main ones are known to be equivalent \cite{gagna2022equivalence}, using $(\infty,2)$-categories as a framework for \emph{formal $\infty$-category theory} is still challenging.

\bigskip
One solution to this problem is represented by Riehl and Verity's theory of $\infty$-cosmoi \cite{RVbook}. In their work, the authors work with a model of $(\infty,2)$-category with enough structure to develop a great amount of the theory of $\infty$-categories in it. Moreover, quite surprisingly, they are able to recover some of this theory by working inside a low-dimensional quotient of the full $\infty$-cosmos, namely its homotopy 2-category. The key point is to study diagrams inside the homotopy 2-category of $\infty$-cosmos (that one might think as \emph{homotopy commutative diagrams}) and then internalise them to the homotopy 2-category of diagrams inside the full $\infty$-cosmos (i.e.\ \emph{homotopy coherent} ones) through a special kind of 2-functor that allows us to lift all the relevant information. This idea is strongly reminiscent of the classical theory of derivators, developed by Grothendieck and Heller \cite{grothendieck1983pursuing, grothendieck1990derivateurs, heller1988homotopy} to describe homotopy limits and colimits in a universal way. In that context, one is interested in capturing what happens when we work with collections of homotopy categories of diagram (model) categories. Even though homotopy limits and colimits can be characterised with a universal property also within $\infty$-category theory, using derivators is still meaningful for at least two reasons:
\begin{enumerate}
    \item this approach doesn't rely on any particular model of $\infty$-category, so it clarifies which results in the theory of $\infty$-categories hold for formal reasons and which ones depend essentially on the features of the model that one is working with,
    \item they provide a more ``user-friendly'' approach to $\infty$-category theory, in the sense that many relevant constructions (most notably the ones available inside stable $\infty$-categories) can be studied already at the level of derivators. This approach has been taken for instance in representation theory (see e.g.\ \cite{sava2022infty}).
\end{enumerate}

Furthermore, in the classical theory of derivators one can associate to every good (e.g.\ complete and cocomplete) $\infty$-category a derivator by looking at the homotopy category of diagrams inside it. With a similar procedure we can associate to every good (e.g.\ combinatorial) model category a derivator. Recalling that all the main models of $\infty$-categories are fibrant objects inside (Quillen-equivalent) model categories (see Figure 1 in the Introduction of \cite{barwick2011unicity}) and that we can associate an $\infty$-category to a model category by simplicially localise it, we have the following picture:
\[\begin{tikzcd}
	{(\infty,1)\text{-categories}} && {\text{Model categories}} \\
	\\
	& {\text{Derivators}}
	\arrow[squiggly, tail reversed, from=1-1, to=1-3]
	\arrow[squiggly, from=1-1, to=3-2]
	\arrow[squiggly, from=1-3, to=3-2]
\end{tikzcd}\]

The natural next chapter in this story is to investigate what a higher version of a derivator would look like. This should be able to blend the axiomatic approach to homotopy theory provided by derivators with the ideas about internalisation that are available in the theory of $\infty$-cosmoi. Hence, it should fit in the following picture:
\[\begin{tikzcd}
	{(\infty,2)\text{-categories}} && {\text{Enriched model categories}} \\
	\\
	& {\text{2-Derivators}}
	\arrow[squiggly, tail reversed, from=1-1, to=1-3]
	\arrow[squiggly, from=1-1, to=3-2]
	\arrow[squiggly, from=1-3, to=3-2]
\end{tikzcd}\]
Using the homotopy hypothesis, which identifies homotopy types with $(\infty,0)$-categories (a.k.a.\ $\infty$-groupoids)\footnote{Caveat: it is really important to be careful about what model for $\infty$-groupoid one is using. See \cite{simpson1998homotopy} for a famous counterexample to one of the models.}, our slogan becomes:
\begin{center}
    2-derivators : $(\infty,1)$-category theory \hspace{5mm} = \hspace{5mm} derivators : homotopy theory
\end{center}
This is justified by the fact that derivators essentially deal with categories, where we can understand ``sameness'' between objects only by means of isomorphisms. In the theory of 2-derivators we are working with 2-categories, so we have much more information available. For instance we can talk about equivalences, adjunctions and so on. In some sense, $2$-derivators provide a way to approach $(\infty,2)$-category theory in a model independent way similarly to how derivators do for $\infty$-categories. Furthermore they are a step forward from the theory of $\infty$-cosmoi in the sense that we never choose a model of $(\infty,2)$-category. To recover the results available in $\infty$-cosmology it is sufficient that an $(\infty,2)$-category is well-behaved enough to give rise to a 2-derivator. 

Fixing a 2-category of categories $\Dia$ to serve as diagram shapes, a derivator is defined as a 2-functor 
$\Dia^{\op}\to\CAT$ satisfying certain axioms. At level 2, one is naturally tempted to define a 2-derivator as a 3-functor from a 3-category of diagram shapes (contained within the 3-category of 2-categories) into the 3-category of 2-categories $\twoCAT$. This choice is debatable, as we explain below.
\subsection{Codomain for higher derivators}\label{intro: codomain}
A derivator is a 2-functor valued in $\CAT$, the (large) 2-category of locally small categories, as this is where the homotopy categories of diagram categories naturally reside. As mentioned earlier, a natural choice for the codomain when studying homotopy 2-categories of enriched diagram categories is $\twoCAT$, the 3-category of 2-categories, 2-functors, 2-natural transformations, and modifications. However, in this work, we instead use $\GRAY$, the (large) Gray-category of 2-categories whose 1-cells are 2-functors, 2-cells are pseudonatural transformations, and 3-cells are modifications. This decision is motivated by two main considerations: first, the examples we examine often involve weak structures and fibrant or cofibrant replacements; second, a 2-natural transformation that is a componentwise equivalence is not generally invertible as a 2-natural transformation, but it is invertible as a pseudonatural transformation. Since in higher derivator theory we would like to be able to characterize equivalences componentwise, just as in ordinary derivators, $\GRAY$ provides a more suitable codomain. Another point that deserves to be made is that we do not require morphisms in the codomain to be weak, even though we work with a weaker notion of adjunction $-$ that is, with biadjunctions. While it is true that biadjoints of 2-functors need not be strict, in our examples these biadjoints arise by enriching ordinary adjunctions. As a result, they are strict 2-functors, and we apply constructions that preserve the strictness of the 2-functors involved. However, these constructions do not preserve strictness at the level of 2-natural transformations.

\subsection{Domain for higher derivators}\label{intro: domain}
A preliminary account of 2-derivators was given in the author’s MRes thesis \cite{NDVmres}, where the 3-category of diagrams is considered as a sub-3-category of 
$\twoCat$, the 3-category of small 2-categories. This choice is inspired by the literature on derivators such as \cite{maltsiniotis2007k}, where the 2-category of diagrams one considers is a full sub-2-category of $\Cat$. In the following we propose a more general higher category of diagrams, in which we consider all simplicially enriched categories. We do so to be able to deal with cotensors with simplicial sets instead of categories. This choice is motivated for instance by 
Riehl and Verity's paper \cite[\S 5]{riehl2015completeness}. In that context the authors describe the forgetful functor from the $\infty$-category of algebras for a homotopy coherent monad by considering cotensors with boundaries of standard $n$-simplices: these are simplicial sets that do not lie in the essential image of the nerve functor $N\colon\Cat\to\sSet$. We conjecture that, in this enlarged diagram category, one can recover the forgetful functor in a similar way.

\section{Background}\label{chap:background}
\graphicspath{{ch_other1/figures/}} 

\subsection{Categorical homotopy theory}
Categorical homotopy theory combines aspects of Homotopy Theory and Category Theory, involving ``categorical approaches to Homotopy Theory'' and ``Homotopy Theory-inspired Category Theory''. In this section, we focus on three intertwined areas: model categories, $\infty$-categories (via quasi-categories), and derivators. For a comprehensive introduction, see Emily Riehl’s book \cite{riehl2014categorical}. Notations may vary, with ordinals denoted as $\mathbb{0}, \mathbb{1}, \mathbb{2}, \dots$ instead of $[-1], [0], [1]$, and so on.

\subsubsection{Model categories}

Model categories appear for the first time in Daniel Quillen's book \emph{Homotopical algebra} \cite{quillen2006homotopical}, where they are called \emph{closed model categories}. They provide a setting for homotopy theory which generalizes many of the classical features of the homotopy theory of (good) topological spaces. Modern textbooks on model category theory include \cite{hirschhorn2003model} and \cite{hovey2007model}. 
	
	\begin{defn}
		A \emph{model category} is a complete and cocomplete category endowed with a model structure, i.e.\ a choice of a triple of subcategories ($\mathcal{W}$, Cof, Fib) of $\M^{\b2}$ resp.\ called \emph{weak equivalences} (w.e.\ for short), \emph{cofibrations} and \emph{fibrations}, and two functorial factorizations $(\alpha,\beta)$ and $(\gamma,\delta)$ satisfying the axioms in Definition 1.1.3. of \cite{hovey2007model}.
		
	\end{defn}
	It is standard practice to denote weak equivalences with $\xrightarrow{\sim}$, fibrations with $\twoheadrightarrow$ and cofibrations with $\hookrightarrow$. 
The appropriate notions of functor and equivalence between model categories are the Quillen ones (see \cite[\S 1]{hovey2007model}).
One of the main advantages of working with model categories is that the process of localizing with respect to weak equivalences has an easy description in terms of homotopies (in the sense of \cite[Definition 1.2.4]{hovey2007model}) between morphisms: the homotopy category $\Ho \M$ of a model category $\M$ having as objects the fibrant-cofibrant objects of $\M$ and as morphisms the homotopy classes of morphisms in $\M$ satisfies the same universal property of $\M[\W^{-1}]$ so the two can be identified. Furthermore Quillen adjunctions and Quillen equivalences descend to adjunctions and equivalences of homotopy categories.
\subsubsection{\texorpdfstring{$\infty$-}{∞-}categories}
As is customary, we call $\infty$-categories what would be more correctly called $(\infty,1)$-categories, so higher categories with $n$-cells (``$n$-dimensional morphisms'') for every $n$ that are invertible in a suitable sense for $n>1$. The literature on $\infty$-categories is quite vast. In this work we will model them using \emph{quasi-categories}, introduced by Boardman and Vogt with the name of \emph{weak Kan complexes} in \cite{boardman2006homotopy}, whose theory has been extensively developed by Joyal and Lurie (see \cite{joyal2008theory}, \cite{lurie2009higher} and \cite{luriehigher}).

Recall that $\Delta$ is the full subcategory of $\Cat$ spanned by the finite non-empty ordinals and $\sSet$ is the category of simplicial sets, namely the category of presheaves of sets over $\Delta$. One of the first results in the theory of simplicial sets establishes the existence of a fully faithful functor $N\colon\Cat\hookrightarrow\sSet$, called the \emph{nerve} (see Chapter 1 of \cite{lurie2009higher}). This functor is part of an adjunction $h\dashv N$, in which the left adjoint is called the \emph{fundamental category} (see \cite{joyal2008theory}) or the homotopy category.

The image of $\Cat$ under the nerve functor sits inside the full subcategory $\qCat\subset\sSet$ spanned by the quasi-categories, in particular it is the subcategory of those quasi-categories that admit \emph{unique} extensions for inner horns (see Chapter 1 of \cite{lurie2009higher}). The following result,  due to Andr\'e Joyal, is also proved by Lurie in \cite{lurie2009higher}.
\begin{thm}
   The category of simplicial sets admits a model structure, called the \emph{Joyal model structure}, in which quasi-categories are the fibrant objects and monomorphisms are the cofibrations.  
\end{thm}

\subsubsection{Derivators}
Derivators were introduced by Grothendieck in \cite{grothendieck1983pursuing} (and further studied in \cite{grothendieck1990derivateurs}). Heller independently developed their theory in his book \cite{heller1988homotopy}, where he calls them \emph{homotopy theories}. Modern references include \cite{groth2013derivators}, \cite{maltsiniotis2001introduction} and \cite{maltsiniotis2007k}.

Let $\Dia$ be a full sub-2-category of $\Cat$, that we usually assume to be closed under some limits and colimits (see \cite{maltsiniotis2007k}). This 2-category can be thought as a 2-category of diagrams, hence the name. In the recent literature (such as \cite{groth2013derivators}) it is common to take $\Dia$ to be the whole of $\Cat$, but we prefer to work in a greater generality.  
\begin{defn}
    A \emph{prederivator} is a 2-functor $\bbD\colon\Dia^{\op}\to\CAT$.
\end{defn}
The image $\bbD(\bb1)$ of the terminal category $\bb1$ is called the \emph{base} of $\bbD$ or \emph{the fiber over $\bb1$}, $\bbD(I)$ is the category of \emph{coherent diagrams}
of shape $I$ or the fiber over $I$ and $\bbD(\bb1)^I$ is the category of \emph{incoherent diagrams} of shape $I$.
\begin{defn}
Given a category $\C$ in $\Dia$, the \emph{represented prederivator} is the 2-functor $\Dia(-,\C)$ .
\end{defn}
\begin{exmp}\label{pder associated to model cat}
If $\M$ is a combinatorial model category and $J$ is a category, the functor category $[J,\M]$ admits both the projective and the injective model structure. Either way, weak equivalences are defined componentwise. The assignment
\begin{align*}
    \Dia^{\op}&\to\CAT\\ 
    J&\mapsto\Ho([J,\M])
\end{align*}
defines a prederivator, called \emph{the prederivator associated to the model category} $\M$. 
\end{exmp}
    
\begin{exmp}
As a consequence of \cite[Remark 15.2.2]{riehl2014categorical}, the simplicial mapping space $\sSet(\Delta^{-} \times \A, \C)$ is a quasi-category whenever $\C$ is so. Hence the assignment
\begin{align*}
    \Dia^{\op}&\to\CAT\\ 
    J&\mapsto h([NJ,\C])
\end{align*}
defines a prederivator\footnote{$h$ is the homotopy category functor introduced in the section on $\infty$-categories.}, called \emph{the prederivator associated to the quasi-category} $\C$. The functoriality of this construction comes from \cite[Theorem 5.14]{joyal2008theory}. 
\end{exmp}
\begin{defn} The diagram 
\begin{center}
\begin{tikzcd}
A \arrow[r, "f^*", ""{name=A, below}] \arrow[d, "g^*"'] & B \arrow[d, "k^*"] 
\ar[ld,shorten <>=10pt,Rightarrow, "\alpha"] 
\\
C \arrow[r, "h^*"', ""{name=C, below}]                  & D 
\end{tikzcd}
\end{center}
in which there are the adjoint pairs $f_!\dashv f^*$ and $h_!\dashv h^*$ is said to be \emph{left Beck-Chevalley} if the mate 
\begin{center}
\begin{tikzcd}
                       & B \ar[rd,shorten <>=10pt,Rightarrow, "\alpha"]\arrow[r, "k^*" ]           & D \ar[rd, xshift=-.5em, yshift=.5em, shorten <>=15pt,Rightarrow, "\epsilon"] \arrow[r, "h_!"]            & C \\
B \arrow[r, "f_!"'] \arrow[ru, bend left, equal, ""{name=A}] & A  \arrow[Rightarrow, shorten <>=7pt, from=A, swap, "\eta"]\arrow[u, "f^*"] \arrow[r, "g^*"'] & C \arrow[u, "h^*"'] \arrow[ru, bend right, equal, ""{name=B}] &  {}
\end{tikzcd}
\end{center}
is an isomorphism. There exists an analogous right Beck-Chevalley condition, which is dual to the one above.
\end{defn}
\begin{rmk}
Recall that whenever $\C$ is a complete and cocomplete category, both left and right Kan extensions of $X\in[J,\C]$ along $u\colon J\to I$ exist and can be computed in a pointwise fashion as $$(u_!X)(k)\coloneqq k^*u_!X\cong\colim_{u/k}\text{pr}^*X \hspace{1cm} \text{and}\hspace{1cm} (u_*X)(k)\coloneqq k^*u_*X\cong\lim_{k/u}\text{pr}^*X,$$ where $k\colon\bb1\to I$ is an object of $I$ and the two ``pr'' stand for the suitable forgetful functors from the comma categories. We can rephrase these formulas in terms of a Beck-Chevalley condition. For instance, the left Kan extension formula is equivalent to saying that
\begin{center}
\begin{tikzcd}
\C^I \arrow[r, "u^*", ""{name=A, below}] \arrow[d, "k^*"'] & \C^J \arrow[d, "\text{pr}^*"] 
\ar[ld,shorten <>=15pt, Rightarrow]\\
\C \arrow[r, "\Delta"', ""{name=C, below}]                  & \C^{(u/k)} 
\end{tikzcd}
\end{center}
is left Beck-Chevalley.
\end{rmk}
\begin{defn}\label{shifted prederivator}
The \emph{shifted prederivator} $\bbD^{I}$ of a prederivator $\bbD$ with respect to $I\in\Dia$
\begin{align*}
    \bbD^{I} \colon\Dia^{\op}&\to\CAT\\ J&\mapsto\bbD(I\times J)
\end{align*}
is obtained as the composition of the $2$-functors $\bbD$ and $I\times -$.
\end{defn}
The shifted prederivator has base $\bbD^I(\bb1)=\bbD(I\times\bb1)\cong\bbD(I)$, which allows us to infer information about any $\bbD(I)$ by working only in the fiber over $\bb1$ and then applying a shift. In fact, as we will see later, the axioms that define a derivator are preserved under the shift operation.

The theory of derivators allows us to recover information about coherent diagrams from what we know at the level of incoherent diagrams. This is achieved through the use of underlying diagram functors.
\begin{defn}\label{underlying diagram functor}
Given $I\in\Dia$, the \emph{underlying diagram functor} 
\begin{align*}
   \dia_I\colon\bbD(I)&\to\bbD(\bb1)^I\\ X&\mapsto\dia_IX
\end{align*}
assigns to every coherent diagram an incoherent one defined as follows:
\begin{align*}
   \dia_IX\colon I&\longrightarrow\bbD(\bb1)\\ i&\longmapsto X_i\coloneqq i^*X\\
   (i\xrightarrow{\alpha}&j)\mapsto (X_i\xrightarrow{X_\alpha} X_j)\coloneqq i^*X\xrightarrow{\alpha^*_X} j^*X
\end{align*}
where we see an object $i\in I$ as a functor $\bb1\xrightarrow{i}I$ and a morphism $i\to j$ as a natural transformation between the respective functors.
\end{defn}
\begin{rmk} In other words, the underlying diagram functors are the ones coming from the action of $\bbD$ on morphisms and the product-internal hom adjunction:
\begin{prooftree}
\AxiomC{$I \xrightarrow{\bbD_{\bb1,I}}\bbD(\bb1)^{\bbD(I)}$}
\UnaryInfC{$\bbD(I)\times I \to\bbD(\bb1)$}
\UnaryInfC{$\dia_I\colon\bbD(I)\to\bbD(\bb1)^{I}$}
\end{prooftree}
\end{rmk}
\begin{defn}
A prederivator $\bbD\colon\Dia^{\op}\to\CAT$ is called a \emph{derivator} if the following axioms hold.

\begin{description}
\item[(Der 1)] $\bbD(\varnothing)=\bb1$ and $\bbD(\coprod_{a\in A}I_a)\xrightarrow{\sim}\prod_{a\in A}\bbD(I_a)$ is an equivalence of categories,
\item[(Der 2)] $\dia_I\colon\bbD(I)\to\bbD(\bb1)^I$ is conservative for every $I$,
\item[(Der 3)] the image $u^*$ through $\bbD$ of every $u\colon J\to K$ in $\Dia$ has both left and right adjoints, called \emph{homotopy left and right Kan extensions} along $u$,
\item[(Der 4)] homotopy Kan extensions are pointwise, i.e.\ the images under $\bbD$ of 
\begin{center}
    \begin{tikzcd}
(u/k) \arrow[r, "\text{pr}"] \arrow[d, "\text{pt}"'] & J\ar[ld,shorten <>=15pt,Rightarrow] \arrow[d, "u"] & (k/u) \arrow[r, "\text{pt}"] \arrow[d, "\text{pr}"'] & \bb1 \ar[ld,shorten <>=15pt,Rightarrow]\arrow[d, "k"] \\
\bb1 \arrow[r, "k"']                      & K                & J \arrow[r, "u"']                      & K               
\end{tikzcd}
\end{center}
are (left and right) Beck-Chevalley squares, where  ``pt'' is the unique functor to the terminal category.
\end{description}
\end{defn}

\begin{prop}
The represented prederivator $\Dia(-,\C)$ is a derivator for every complete and cocomplete category $\C$.
\end{prop}

\begin{thm}\label{shider is der} Let $\bbD$ be a derivator and $I\in\Dia$. Then the shifted prederivator $\bbD^I\colon\Dia^{\op}\to\CAT$ is a derivator.
\end{thm}

An interesting subclass of derivators satisfies another axiom. Such derivators are known as \emph{strong derivators}.
\begin{axiom}[(Der 5)]\label{Der 5}
 Let $\b2$ be the walking arrow. The functor $\dia_{\b2}\colon\bbD^I(\b2)\to(\bbD^I(\bb1))^{\b2}$ is full and essentially surjective for every $I\in\Dia$.
\end{axiom}
\begin{rmk}
Since $\bbD^I$ is a derivator whenever $\bbD$ is, the functor involved in \ref{Der 5} is conservative. That is, strong derivators are the derivators for which the functor  $\dia_{\b2}\colon\bbD^I(\b2)\to(\bbD^I(\bb1))^{\b2}$ is \emph{weakly smothering} for every $I\in\Dia$.
\end{rmk}
Unsurprisingly this good property is enjoyed by a large class of derivators, in particular by those associated to nice model categories as showed in \cite{groth2013derivators}. 
\begin{prop}
Let $\M$ be a combinatorial model category. Then the prederivator associated to $\M$ (see Example \ref{pder associated to model cat})
is a strong derivator.
\end{prop}
\subsection{Change of base}
If we are given a $\V$-category $\C$ and a monoidal functor $T\colon\V\to\W$, there is a way to use $T$ to change base of enrichment for $\C$. This is what is called \emph{change of base}. We briefly recall a few basic facts, referring to \cite[Ch.\ 4]{cruttwell2008normed}, \cite[\S II]{EilenbergS:cloc} and \cite[A.6.]{more_elements} for proofs.

\begin{prop}
     A lax monoidal functor $T\colon\V\to\W$ induces a change-of-base 2-functor
    $$T_*\colon\VCat\to\WCat$$
   defined on objects as
    \begin{enumerate}
       \item $\Ob(T_*\C)=\Ob(\C)$ for every $\V$-category $\C$,
        \item $(T_*\C)(X,Y)=T(\C(X,Y))$, with compositions and identity being the composites of the morphisms arising from the monoidality of $T$, as well as the composition and unit of the 
$\V$-category $\C$,
\[\begin{tikzcd}[column sep= .8 em]
	{T\C(Y,Z)\otimes_{\W} T\C(X,Y)} && {T(\C(Y,Z)\otimes_{\V}\C(X,Y))} && {T\C(X,Z)}
	\arrow[from=1-1, to=1-3]
	\arrow[from=1-3, to=1-5]
\end{tikzcd}\]
\[\begin{tikzcd}[column sep= .8 em]
	{I_{\W}} && {TI_{\V}} && {T\C(X,X),}
	\arrow[from=1-1, to=1-3]
	\arrow[from=1-3, to=1-5]
\end{tikzcd}\]
    \end{enumerate}
sending a $\V$-functor $F\colon\C\to\D$ to the $\W$-functor $T_*F\colon T_*\C\to T_*\D$ acting as $F$ on objects and as
$$T\C(X,Y)\xrightarrow{TF}\D(TX,TY)$$
on hom-objects. Finally, a $\V$-natural transformation $\alpha\colon F\Rightarrow G$ is sent to a $\W$-natural transformation $T_*\alpha\colon T_*F\Rightarrow T_*G$ with components
$$ I_{\W}\rightarrow TI_{\V}\rightarrow T\D(Fc,Gc).$$
coming from the monoidality of $T$ and the components of $\alpha$.
\end{prop}
\begin{exmp}
    For every cartesian closed category $(\V,\times, I)$, the functor $(-)_0\coloneqq\V(I,-)\colon\V\to\Set$ is strong monoidal (since covariant homs preserve limits) hence inducing a change of base 2-functor $\VCat\to\Cat$ giving the \emph{underlying category} of a $\V$-category.
\end{exmp}
\begin{exmp}\label{homotopy 2-cat}
The homotopy category functor $h\colon\sSet\to\Cat$ preserves products therefore it induces a 2-functor $h_*\colon\sSet\mbox{-}\Cat\to\twoCat$ which returns the \emph{homotopy 2-category} of a simplicially enriched category. The homotopy 2-category of a $\sSet$-category $\M$ has the same objects of $\M$, $0$-simplices in the hom-simplicial sets of $\M$ as $1$-cells and homotopy classes of formal
composites of $1$-simplices in the hom-simplicial sets of $\M$ as 2-cells. When we restrict to $\qCat$-categories, the $2$-cells are just the homotopy classes of $1$-simplices in the hom-simplicial sets of $\M$.
\end{exmp}
As a Corollary of the 2-dimensional aspects of sending a base monoidal category $\V$ to $\V\mbox{-}\Cat$, we recover the following Proposition.
\begin{prop}
    Given cartesian closed categories $\V,\W$ and an adjunction $F\dashv G$ between them, where both $F$ and $G$ preserve finite products, we get a change-of-base 2-adjunction $F_*\dashv G_*$.
\end{prop}
For instance, the adjunction $h\dashv N$ between the homotopy category and the nerve functor induce a change of base 2-adjunction $h_*\colon\sSet\mbox{-}\Cat\rightleftarrows\twoCat\colon N_*$. Since $N_*$ is also fully faithful, in the following we will routinely identify a 2-category $\C$ with its simplicially enriched
counterpart $N_*\C$, dropping the extra $N_*$ to simplify the notation.
\subsection{Weighted limits via collages}

All the important concepts from category theory have enriched counterparts. In the case of limits and colimits, the appropriate notions are \emph{weighted limits} and \emph{weighted colimits}. We refer to Section 3.9 of Max Kelly's book \cite{kelly1982basic} for a discussion of these and an explanation of why they are the correct notions. Among weighted limits and colimits, we recall that pointwise Kan extensions in the $2$-category $\VCat$, when they exist, can be computed as such. We follow Kelly's book in working exclusively with pointwise Kan extensions, and therefore drop the adjective ``pointwise" in the following definition.
\begin{defn}
The \emph{left Kan extension} of a $\V$-functor $G\colon\A\to\B$ along a $\V$-functor $K\colon\A\to\C$, if it exists, is a $\V$-functor $\Lan_{K}G\colon\C\to\B$ defined by the formula
$$\Lan_{K}G(c)=\colim(\C(K-,c),G)$$
Dually, the \emph{right Kan extension} is defined by the formula
$$\Ran_{K}G(c)=\lim(\C(c,K-),G).$$
\end{defn}

Another way to express weighted limits is via their \emph{collage}, which makes sense for a more general \emph{enriched profunctor}. 
\begin{defn}
An \emph{enriched profunctor} or \emph{$\V$-profunctor} $W\colon\A\xslashedrightarrow{}\B$ between the $\V$-categories $\A$ and $\B$ is a $\V$-functor $\B^{\op}\otimes\A\to\V$.
\end{defn}
\begin{exmp}
    $\Set$-profunctors are the usual profunctors between locally small categories $\A$ and $\B$, namely functors $\B^{\op}\times\A\to\Set$. 
\end{exmp}
      
\begin{rmk}\label{enriched homs}
    Using the properties of the unit $\V$-category $\I$, a weight $\A\to\V$ is the same as a $\V$-functor $\I^{\op}\otimes\A\to\V$ or in other words a $\V$-profunctor $\A\xslashedrightarrow{}\I$.
\end{rmk}




\begin{defn}\label{collage}
Assume that $\V$ has an initial object $\varnothing$ that is preserved on both sides by $\otimes$. Given a $\V$-profunctor $W\colon\A\xslashedrightarrow{}\B$, we define its \emph{collage} $\coll(W)$ to be the $\V$-category having as objects the coproduct $\Ob\A\sqcup\Ob\B$ and such that 
    \[
\coll(W)(x,y)=
\begin{cases} 
\A(x,y), & \mbox{if } x,y\in\A \\ 
\B(x,y), & \mbox{if } x,y\in\B\\ 
W(x,y), & \mbox{if } x\in\B \ \text{and}\ y\in\A\\ 
\varnothing, & \mbox{otherwise }
\end{cases} 
\]
where the composition comes from the ones in $\A$ and $\B$ and from the $\V$-functoriality of $W$.
\end{defn}
As a special case, we recover the notion of collage of a weight, described e.g.\ in \cite[Definition 7.2.7]{more_elements} for $\V=\sSet$.

From the definition we obtain inclusions $i_{\A}^{W}\colon\A\hookrightarrow\coll(W)$ and $i_{\B}^{W}\colon\B\hookrightarrow\coll(W)$, forming a cospan $\A\hookrightarrow\coll(W)\hookleftarrow\B$ in the $2$-category of $\V$-categories. In most cases, cospans coming from collages of profunctors are exactly the \emph{two-sided codiscrete cofibrations} (in the sense of \cite{street1980fibrations}) in the $2$-category in which they live. This inspires an alternative approach to the collage construction using cospans.  
\begin{defn}\label{collage as cospan}
We define the \emph{collage} $\coll(f,g)$ of the cospan $\A\xrightarrow{f}\C\xleftarrow{g}\B$ in $\V\mbox{-}\Cat$ as the $\V$-category with set of objects the coproduct $\Ob\A\sqcup\Ob\B$, hom-objects $\A(a,a')$ and $\B(b,b')$ between $a,a'\in\A$ and $b,b'\in\B$, $\C(fa,gb)$ from an object $a\in\A$ to an object $b\in\B$ and $\varnothing$ elsewhere.
\end{defn}
 
Considering the $\V$-profunctor $\C(f-,g-)\colon\B\xslashedrightarrow{}\A$ we can recover Definition \ref{collage as cospan} from Definition \ref{collage}. In general, there is a $\V$-functor $\pi_{\C}\colon\coll(f,g)\to\C$ acting on objects as the coproduct $f\sqcup g$, on hom-objects as $f$ in the full subcategory $\A$, as $g$ in the full subcategory $\B$, as the identity on hom-objects from an element of $\A$ to an element of $\B$ and as the unique  morphism $\varnothing=\coll(f,g)(b,a)\to\C(\pi_{\C}b,\pi_{\C}a)=\C(gb,fa)$ from $b\in\B$ to $a\in\A$. 

\begin{rmk}\label{collage is idempotent}
    The process of taking collages is idempotent, in the sense that the collage of the cospan $\A\xrightarrow{f}\C\xleftarrow{g}\B$ in $\V\mbox{-}\Cat$ gives a new cospan $\A\xhookrightarrow{i_{\A}}\coll(f,g)\xhookleftarrow{i_{\B}}\B$ in $\V\mbox{-}\Cat$ whose collage is again $\coll(f,g)$. Notice indeed that $\coll(i_{\A},i_{\B})=\coll(f,g)$ since they have the same set of objects and 
\begin{align*}
    \coll(i_{\A},i_{\B})(a,b)&=\coll(f,g)(i_{\A}a,i_{\B}b) \\
    &=\coll(f,g)(a,b)\\
    &=\C(fa,gb)\\
\end{align*}
\end{rmk}
We conclude this section by explaining how to compute weighted colimits using collages. This is achieved by the the following characterization, a proof of which for $\V=\sSet$ can be found in \cite[\S 7.2]{more_elements}.
\begin{prop}\label{weighted limits as Kan extensions}
Let $F\colon\A\to\B$ be a $\V$-functor and $W\colon\A\to\V$ a weight. The weighted limit $\lim(W,F)$ exists if and only if the pointwise right Kan extension of $F$ along the inclusion $i_{\A}^{W}\colon\A\hookrightarrow\emph{coll}(W)$ exists. In this case it can be computed as $$\lim(W,F)\cong\Ran_{i_{\A}^{W}}(\bullet),$$ where $\bullet$ is the unique object of $\bb1$ in the collage. Dually, we can compute weighted colimits as \emph{left} Kan extensions.
\end{prop}
\subsection{A brief excursion in enriched model category theory}
As outlined in the Introduction, transitioning from derivators to 2-derivators corresponds to moving from diagrams within a model category to diagrams in an enriched model category. Enriched model categories are, in general, model categories that allow for an enrichment that is compatible with the model structure (see \cite{hovey2007model, lurie2009higher} for precise definitions). The reference text for enriched category theory is Kelly's book \cite{kelly1982basic}. 
\begin{exmp}
The Quillen model structure on simplicial sets is a closed symmetric monoidal model category. Model categories enriched over it are usually called \emph{simplicial model categories}. 
\end{exmp}
\begin{exmp}
The Joyal model structure is a monoidal model category with the monoidal product being the cartesian product. It is then self-enriched as a model category, meaning that it is $\sSet_{\joy}$-enriched and thus presents an $(\infty,2)$-category\footnote{This is analogous to the definition of a $2$-category as a $\Cat$-enriched category.}.
\end{exmp}

\begin{rmk}\label{weighted limits as adjoints}
Let $\C,\D$ be $\V$-categories and $W\colon\D\to\V$ a weight. Suppose moreover that $\C$ has tensors. We can then characterize a $W$-weighted limit in $\C$ as the right adjoint to the functor 
\begin{align*}
    W\otimes -\colon\C&\to[\D,\C]\\
    c&\mapsto W\otimes c
\end{align*}
where
\begin{align*}
    W\otimes c\colon\D&\to\C\\
    d&\mapsto Wd\otimes c.
\end{align*}
The right adjoint to $W\otimes-$ sends indeed a functor $F\in[\D,\C]$ to an object $x\in\C$ such that $[\D,\C](W\otimes c, F)\cong\C(c,x)$. But since $\C$ is tensored (and then so is $[\D,\C]$), we have the isomorphism $[\D,\C](W\otimes c, F)\cong[\D,\V](W,\C(x,F))$ that implies $x\cong\lim^{W}F$.
\end{rmk}
The following result is folklore, and we just adapted it to our purposes (for detailed explanations of the terminology used, see \cite{moser2019injective}).
\begin{prop}\label{inj equiv proj} Let $\M$ be a locally $\V$-presentable $\V$-category admitting an accessible $\V$-enriched model structure and take $\D$ to be a small $\V$-category. Suppose that $\D(d,d')$ is a cofibrant object of $\V$ for every couple of objects $d,d'\in\D$. Then the adjunction $\id\colon[\D,\M]^{\emph{proj}}\rightleftarrows[\D,\M]^{\emph{inj}}\colon\id$ is a Quillen equivalence.
\end{prop}
\begin{proof}
We have to show that the identity functor $[\D,\M]^\inj\to [\D,\M]^\proj$ is right Quillen, i.e.\ it preserves fibrations and trivial fibrations. That is, we have to prove the following two statements:
\begin{enumerate}[(i)]
    \item every injective fibration\footnote{A morphism that has the RLP with respect to componentwise trivial cofibrations} $\alpha\colon F\to G$ in $[\D,\M]$ is also a projective fibration (i.e.\ a componentwise fibration),
    \item every injective trivial fibration\footnote{A morphism that has the RLP with respect to componentwise cofibrations} $\beta\colon H\to K$ in $[\D,\M]$ is also a projective trivial fibration (i.e.\ a componentwise trivial fibration).
\end{enumerate}
In order to prove (i), let us take an injective fibration $\alpha$. Fix $d\in\D$. We claim that $\alpha_d\colon Fd\to Gd$ has the RLP with respect to every trivial cofibration $i\colon A\to B$ in $\M$. The following square of solid arrows admits a dashed lift
\begin{center}
\begin{tikzcd}
A \arrow[r] \arrow[d, "i"']          & Fd \arrow[d, "\alpha_d"] \\
B \arrow[r] \arrow[ru, dashed] & Gd          
\end{tikzcd}
\end{center}
if and only if
\begin{center}
\begin{tikzcd}
A \arrow[r] \arrow[d, "i"']          & \lim^{\D(d,-)}F\arrow[d, "\widetilde{\alpha}_d"] \\
B \arrow[r] \arrow[ru, dashed] & \lim^{\D(d,-)}G         
\end{tikzcd}
\end{center}
does, since $\lim^{\D(d,-)}F\cong Fd$ by the Yoneda Lemma. Using the adjunction $-\otimes\D(d,-)\dashv\lim^{\D(d,-)}$ from Remark \ref{weighted limits as adjoints}, we can transpose the latter lifting problem to the following one
\begin{center}
\begin{tikzcd}
A \otimes\D(d,-) \arrow[r] \arrow[d, "\operatorname{\emph{i}\hspace{.5mm} \otimes\D(\emph{d}\hspace{.5mm},-)}"']          & F\arrow[d, "\alpha"] \\
B \otimes\D(d,-) \arrow[r] \arrow[ru, dashed] & G.        
\end{tikzcd}
\end{center}
Since $i$ is a trivial cofibration in $\M$ and $\emptyset\to\D(d,d')$ is a cofibration in $\V$ (because by assumption every hom-object is cofibrant), the pushout-product axiom for enriched model categories guarantees that the morphism $$A\otimes\D(d,d')\coprod_{A\otimes\emptyset}B\otimes\emptyset\cong A\otimes\D(d,d')\xrightarrow{i\otimes\D(d,d')} B\otimes\D(d,d')$$is a trivial cofibration in $\M$. This means that $A\otimes\D(d,-)\xrightarrow{i\otimes\D(d,-)} B\otimes\D(d,-)$ is a componentwise trivial cofibration. Thus the last lifting problem can be solved since $\alpha$ has the RLP with respect to componentwise trivial cofibrations, showing (i). The proof of (ii) is exactly the same, with $i\colon A\to B$ a generic cofibration in $\M$.
\end{proof}
As a consequence of Proposition \ref{inj equiv proj} we have the following result.
\begin{cor}\label{equiv proj inj joyal enriched}
Let $\M$ be a $\sSet_{\emph{Joyal}}$-enriched model category. For every small simplicial category $\J$, the projective and the injective enriched model structure on the diagram category $[\J,\M]$, if they exist, are Quillen equivalent. 
\end{cor}
\begin{proof}
This holds since in $\sSet_\joy$ every mono is a cofibration; hence $\varnothing\to\D(d,d')$ is a cofibration for every $d,d'\in\D$. We can then apply Proposition \ref{inj equiv proj} to conclude. 
\end{proof}


\begin{lemma}\label{componentwise equiv} 
Let $\A$ and $\B$ be $\sSet_{\emph{Joyal}}$-enriched model categories, $F, G\colon\A\rightrightarrows\B$ simplicial functors and $\alpha\colon F\Rightarrow G$ a simplicial natural transformation which is a componentwise homotopy equivalence. Then the $2$-natural transformation $h_*\alpha$ is a componentwise equivalence.
\end{lemma}
\begin{proof}
By assumption we have that for every $X\in\A$ there exists $l_X\colon GX\to FX$ s.t.\ $l_X\alpha_X\sim\id_{FX}$ and $\alpha_X l_X\sim\id_{GX}$. In the Joyal model structure the cylinder objects are of the form $\text{Cyl}(Y)=Y\times J$ where $J$ is the nerve of the free-living isomorphism $\bbI$. Hence we can write the first homotopy as
\begin{center}
    \begin{tikzcd}
\Delta^0 \arrow[r, "0"] \arrow[rd, "l_X\alpha_X"'] & J \arrow[d, "H"] & \Delta^0 \arrow[l, "1"'] \arrow[ld, "\id_{FX}"] \\
                                                   & {\B(FX,FX)}      &                                                
\end{tikzcd}
\end{center}
since $\text{Cyl}(\Delta^0)=\Delta^0\times J\cong J$. Furthermore, we have denoted by $0,1$ the two $0$-simplices of $J$. Now we can apply $h$ and we get
\begin{center}
    \begin{tikzcd}
{[0]} \arrow[r, "0"] \arrow[rd, "h(l_X)h(\alpha_X)"'] & \bbI \arrow[d, "hH"] & {[0]} \arrow[l, "1"'] \arrow[ld, "\id_{FX}"] \\
                                                      & {h\B(FX,FX)}      &                                             
\end{tikzcd}
\end{center}
where we used the isomorphism $hN\cong\id_{\Cat}$ and the functoriality of $h$ (notice that $hFX=FX$). In other words we found an invertible $1$-cell inside $h\B(FX,FX)$ whose endpoints are $h(l_X)h(\alpha_X)$ and $\id_{FX}$, meaning that they are isomorphic. With the same argument one shows that the other homotopy induces an isomorphism as well. To conclude the proof we see that $(h_*\alpha)_X\cong h(\alpha_X)$ thanks to the action of the change of base $2$-functor $h_*$ on the monoidal units of the base monoidal categories: it is precomposed by a map $[0]\to h\Delta^0$ which is necessarily an isomorphism since the endpoints are both trivial categories.
\end{proof}
\begin{cor}\label{2-nat ptwise we between fib-cof}
Suppose we have a pair of $2$-functors $G,G'\colon\D\rightrightarrows h_*(\M)$ between a $2$-category $\D$ and a $\sSet_{\emph{Joyal}}$-enriched model category $\M$ such that $GD, G'D$ are fibrant-cofibrant objects in $\M$ for every $D\in\D$. Let $\alpha\colon G\Rightarrow G'$ be a $2$-natural transformation s.t.\ $\alpha_D$ is a weak equivalence for every $D\in\D$. Then $\alpha$ is a pseudonatural equivalence. 
\end{cor}
\begin{proof}
The proof is an  immediate consequence of Lemma \ref{componentwise equiv} and the fact that every 2-natural transformation that is a componentwise equivalence is a pseudonatural equivalence, since the components of $\alpha$ are weak equivalences between fibrant-cofibrant objects.
\end{proof}
\begin{defn}\label{simplicial replacements}
Let $\M$ be a $\sSet_\joy$-enriched model category. We define $F^{\M}, C^{\M}\colon\M\to\M$ to be the fibrant and cofibrant replacement simplicial functors. They come with simplicial natural transformations $\eta^{\M}\from\id_{\M}\Rightarrow F^{\M}$ and $\epsilon^{\M}\from C^{\M}\Rightarrow\id_{\M}$ whose components are trivial cofibrations and trivial fibrations, respectively. When restricted to fibrant-cofibrant objects, they admit pseudonatural inverses. 
\end{defn}

In the following, we generally omit the superscripts from the fibrant and cofibrant replacement simplicial functors to avoid overloading the notation. 

\begin{cor}\label{Cofib and fib replacement composed with pointwise equiv}
Suppose we have a pair of $2$-functors $H,H'\colon\D\rightrightarrows h_*(\M)$ between a $2$-category $\D$ and a $\sSet_{\emph{Joyal}}$-enriched model category $\M$. Let $\alpha\colon H\Rightarrow H'$ be a $2$-natural transformation s.t.\ $\alpha_D$ is a weak equivalence and $HD, H'D$ are fibrant objects in $\M$ for every $D\in\D$. Then $C\alpha$ is a pseudonatural equivalence. Dually, if $HD$ and $H'D$ are cofibrant for every $D\in\D$ then $F\alpha$ is a pseudonatural equivalence.
\end{cor}
\begin{proof}
From the functorial factorizations of the model category $\M$ we know that for every $D\in\D$ the diagram
\[
\begin{tikzcd}
CHD \arrow[r, "\sim", two heads] \arrow[d, "(C\alpha)_D"'] & HD \arrow[d, "\alpha_D"] \\
CH'D \arrow[r, "\sim"', two heads]                         & H'D                      
\end{tikzcd}
\]
commutes and $\alpha_D$ is a weak equivalence. Then $(C\alpha)_D$ is a weak equivalence (for the $2$-out-of-$3$ property of weak equivalences) between fibrant-cofibrant objects (because the cofibrant replacement functor sends fibrant objects to fibrant-cofibrant objects). Therefore Corollary \ref{2-nat ptwise we between fib-cof} implies that $C\alpha$ is a pseudonatural equivalence.  
\end{proof}

An interesting class of enriched functors between enriched model categories is given by enriched Quillen adjunctions, which are enriched adjunctions whose underlying ordinary adjunction is a Quillen adjunction between ordinary model categories. An important fact is that such adjunctions induce biadjunctions between the homotopy 2-categories of fibrant and cofibrant objects, as we will see later. Before that, however, let us recall what a biadjunction is and some facts about them.

\begin{defn}\label{defn biadjunction}
A \emph{biadjunction} $f\dashv_b u$ consists of $2$-categories $\A$ and $\B$, $2$-functors $f\colon\B\to\A$ and $u\colon\A\to\B$ and pseudonatural transformations $\eta\colon\id_{\B}\Rightarrow uf$ and $\epsilon\colon fu\Rightarrow\id_{\A}$ satisfying the triangular identities up to invertible modifications, i.e.\ such that there exist invertible modifications filling the triangles
\[
\begin{tikzcd}
f\arrow[r, Rightarrow, "f\eta"] \arrow[rd, equal] & fuf \ar[ld, phantom, "\hspace{.5cm}\cong", near start]
\arrow[d, Rightarrow, "\epsilon f"] & u \arrow[r, Rightarrow, "\eta u"] \arrow[rd, equal] & ufu \ar[ld, phantom, "\hspace{.5cm}\cong", near start]
\arrow[d, Rightarrow, "u\epsilon"] \\
                         {}             & f                       &  {}                                   & u                    
\end{tikzcd}
\]
\end{defn}
Similarly to ordinary adjunctions, biadjunctions can also be characterized homwise by specifying for each object $a$ in $\A$ and each object $b$ in $\B$ an equivalence of categories $\A(fb,a)\cong\B(b,ua)$, which is pseudonatural both in $a$ and in $b$. These two definitions are equivalent by the Yoneda lemma for bicategories (see e.g.\ Chapter 8 of \cite{JohYau:2-cat}). Next, we provide a criterion for detecting biadjunctions that generalizes a perhaps little-known result from ordinary category theory, which we recall below.
\begin{prop}\label{adj from composites}
Let $g\colon A\to B$ and $f\colon B\to A$ be functors between the categories $A$ and $B$. Suppose there exist two natural transformations $\eta\colon\id_B\Rightarrow gf$ and $\epsilon\colon fg\Rightarrow\id_A$ such that the triangle composites $\epsilon f \cdot f\eta$ and $g\epsilon\cdot\eta g$ are natural isomorphisms. Then $f$ is left adjoint to $g$. 
\end{prop}
\begin{proof}
We need to show that the triangular identities hold, namely we have to prove 
\[
\begin{tikzcd}
f\arrow[r, Rightarrow, "f\eta'"] \arrow[rd, equal] & fgf 
\arrow[d, Rightarrow, "\epsilon' f"] & g \arrow[r, Rightarrow, "\eta' g"] \arrow[rd, equal] & gfg 
\arrow[d, Rightarrow, "g\epsilon'"] \\
                         {}             & f                       &  {}                                   & g                    
\end{tikzcd}
\]
commute for a suitable choice of $\eta'$ and $\epsilon'$. If we call $\Phi\coloneqq\epsilon f \cdot f\eta$ and $\Theta\coloneqq g\epsilon\cdot\eta g$, we have $\Phi^{-1}\Phi=\id_f=\Phi\Phi^{-1}$ and $\Theta^{-1}\Theta=\id_g=\Theta\Theta^{-1}$. Taking $\eta'\coloneqq\Theta^{-1}f\cdot\eta$ and $\epsilon'\coloneqq\epsilon$, we observe that 
\[
\begin{tikzcd}
\id_B \arrow[r, "\eta", Rightarrow] \arrow[d, "\eta "', Rightarrow]                                                     & gf \arrow[r, "g\Phi^{-1}", Rightarrow] \arrow[d, "\eta gf"', Rightarrow] 
& gf \arrow[d, "\eta gf", Rightarrow] \arrow[rdd, "\id_{gf}", Rightarrow, bend left] 
&    \\
gf \arrow[r, "gf\eta", Rightarrow] 
\arrow[rrd, "\id_{gf}"', Rightarrow, bend right=60] & gfgf \arrow[r, "gfg\Phi^{-1}", Rightarrow] \arrow[d, "g\epsilon f"', Rightarrow]                         & gfgf \arrow[d, "g\epsilon f", Rightarrow] 
&    \\
                                                                                                                     & gf\arrow[r, "g\Phi^{-1}"', Rightarrow]                                                                  & gf \arrow[r, "\Theta^{-1}f"', Rightarrow]                                                                  & gf
\end{tikzcd}
\]
is commutative by the middle four interchange, hence $\Theta^{-1}f\cdot\eta=g\Phi^{-1}\cdot\eta$. Using this, we have that the following diagrams
\[
\begin{tikzcd}
                                                                                                                &                                                                                     &                                                                         &                                                                                                                    & {}                                                                                                                  &                                                                        \\
f \arrow[r, "f\eta", Rightarrow] \arrow[rd, "\Phi"', Rightarrow] \arrow[rr, "f\eta'", Rightarrow, bend left=49] & fgf \arrow[d, "\epsilon f", Rightarrow] \arrow[r, "fg\Phi^{-1}", Rightarrow] & fgf \arrow[d, "\epsilon f", Rightarrow] 
& g \arrow[r, "\eta g", Rightarrow] \arrow[rd, "\Theta"', Rightarrow] \arrow[rr, "\eta'g", Rightarrow, bend left=49] & gfg \arrow[d, "g\epsilon", Rightarrow] \arrow[r, "\Theta^{-1}fg", Rightarrow] 
& gfg \arrow[d, "g\epsilon", Rightarrow] 
\\
                                                                                                                & f \arrow[r, "\Phi^{-1}"', Rightarrow]                                        & f                                                                      &                                                                                                                    & g \arrow[r, "\Theta^{-1}"', Rightarrow]                                                                      & g                                                                     
\end{tikzcd}
\]
commute. Since $\Theta^{-1}\Theta=\id_g$ and $\Phi^{-1}\Phi=\id_f$ we get the result.
\end{proof}
\begin{lemma}\label{lemma biadjunction}
Suppose we have a pair of $2$-functors
\[
\begin{tikzcd}
\A\arrow[r, bend right,"u"',""{name=A, below}] & \B\arrow[l, bend right,"f"',""{name=B,above}]
\end{tikzcd}
\]
and a couple of pseudonatural transformations $\eta\colon\id_{\B}\Rightarrow uf$ and $\epsilon\colon fu\Rightarrow\id_{\A}$. If the triangle composites $\epsilon f \cdot f\eta$ and $u\epsilon\cdot\eta u$ are pseudonatural equivalences, then $f$ and $u$ can be promoted to a biadjunction.
\end{lemma}
\begin{proof}
The proof proceeds in the same way as in Proposition \ref{adj from composites}, with the only differences being that we have pseudonatural equivalences instead of natural isomorphisms and diagrams commuting up to invertible modifications instead of strictly commutative ones. 
\end{proof}

\begin{prop}\label{enriched Quillen induces biadj} An enriched Quillen adjunction
\[
\begin{tikzcd}
\M\arrow[r, bend right,"U"',""{name=A, below}] & \cN\arrow[l, bend right,"V"',""{name=B,above}]\ \ar[from=A, to=B, symbol=\vdash]
\end{tikzcd}
\]
between two $\sSet_{\emph{Joyal}}$-enriched model categories induces a biadjunction
\[
\begin{tikzcd}
h_*\M_{\emph{cf}}\arrow[r, bend right,"\ou"',""{name=A, below}] & h_*\cN_{\emph{cf}}\arrow[l, bend right,"\ov"',""{name=B,above}]\ \ar[from=A, to=B, symbol=\perp_{\emph{b}}, rotate=90]
\end{tikzcd}
\]
between the homotopy $2$-categories of fibrant-cofibrant objects.
\end{prop}
\begin{proof}
Since $U$ and $C$ preserve fibrant objects\footnote{In fact, $U$ is right Quillen and $\epsilon^{\cN}_X\from CX\to X$ is a trivial fibration for every $X\in\cN$.}, they can be restricted and composed to give a simplicial functor $$\M_{\text{cf}}\xrightarrow{U}\cN_{\text{f}}\xrightarrow{C}\cN_{\text{cf}}.$$ Similarly, we get a simplicial functor $$\cN_{\text{cf}}\xrightarrow{V}\M_{\text{c}}\xrightarrow{F}\M_{\text{cf}}.$$ Applying $h_*$ we find $2$-functors $$\ou\coloneqq h_*(CU)\from h_*\M_{\text{cf}}\to h_*\cN_{\text{cf}}$$ and $$\ov\coloneqq h_*(FV)\from h_*\cN_{\text{cf}}\to h_*\M_{\text{cf}}.$$ Moreover, the restrictions of the simplicial natural transformations $\eta^{\M}$ 
and $\epsilon^{\cN}$ to the full simplicial subcategories $\M_{\text{cf}}$ and $\cN_{\text{cf}}$, spanned by fibrant-cofibrant objects, are componentwise homotopy equivalences. Therefore, by Lemma \ref{componentwise equiv}, 
$h_*$ sends them to $2$-natural transformations that are componentwise equivalences, and hence admit a (global) pseudonatural inverse. We denote these inverses by $\widetilde{\eta^{\M}}$ and $\widetilde{\epsilon^{\cN}}$ and we define the following pseudonatural transformations 
$$\eta'\coloneqq \id_{h_*{\cN_{cf}}}\xRightarrow{\widetilde{\epsilon^{\cN}}}C\xRightarrow{\overline{\eta}}\ou\ov \ \ \text{and} \ \ \epsilon'\coloneqq \ov\ou \xRightarrow{\overline{\epsilon}}F\xRightarrow{\widetilde{\eta^{\M}}}\id_{h_*{\cM_{cf}}}$$
where $\overline{\eta}$ is the image through $h_*$ of the restriction of the composition $$C\xRightarrow{C\eta}CUV\xRightarrow{CU\eta^{\M} V} CUFV$$ to $\cN_{cf}$ and $\overline{\epsilon}$ is obtained from the composition $$FVCU\xRightarrow{FV\epsilon^{\cN} U}FVU\xRightarrow{F\epsilon}F$$ in a similar way ($\eta$ and $\epsilon$ are the unit and the counit of the enriched Quillen adjunction $V\dashv U$). It remains to show that $\epsilon'\ov \cdot \ov\eta'$ and $\ou\epsilon'\cdot \eta'\ou$ are pseudonatural equivalences, so that we can apply Lemma \ref{lemma biadjunction} to get the claim. First notice that the diagram
\begin{center}
    \begin{tikzcd}[
  column sep={10em,between origins},
  row sep={5em,between origins},
]
CCU \arrow[r, "C \eta CU", Rightarrow] \arrow[d, "\epsilon^{\cN} CU"', Rightarrow]                   & CUVCU \arrow[d, "CUV\epsilon^{\cN} U" description, Rightarrow] \arrow[r, "CU\eta^{\M} VCU", Rightarrow]  & CUFVCU \arrow[d, "CUFV\epsilon^{\cN}  U", Rightarrow] \\
CU \arrow[r, "C\eta U"', Rightarrow] \arrow[rd, equal, bend right] & CUVU \arrow[r, "CU\eta^{\M} VU"', Rightarrow] \arrow[d, "CU\epsilon" description, Rightarrow]                                  & CUFVU \arrow[d, "CUF\epsilon", Rightarrow]   \\
                                                                                            & CU \arrow[r, "CU \eta^{\M}"', Rightarrow]                                                                                      & CUF                                                                        
\end{tikzcd}
\end{center}
is made of squares which are commutative thanks to the middle four interchange (they are $2$-cells in the $2$-category of simplicial categories, simplicial functors and simplicial natural transformations) and a triangle that is commutative by one of the triangular identities of the adjunction $V\dashv U$. Using the middle four interchange we find that the diagram
\begin{center}
    \begin{tikzcd}
CCU \arrow[r, "\epsilon^{\cN} CU", Rightarrow] \arrow[d, "\epsilon^{\cN} CU"', Rightarrow] & CU \arrow[r, "CU\eta^{\M}", Rightarrow] \arrow[d, "\epsilon^{\cN} U", Rightarrow] & CUF \arrow[d, "\epsilon^{\cN} UF", Rightarrow] \\
CU \arrow[r, "\epsilon^{\cN} U"', Rightarrow]                                     & U \arrow[r, "U\eta^{\M}"', Rightarrow]                                   & UF                                   
\end{tikzcd}
\end{center}
is commutative as well. The vertical $2$-cells and the horizontal $2$-cells of left square are componentwise weak equivalences because $\epsilon^{\cN}$ is so. In addition, $U$ preserves trivial fibrations (it is right Quillen) and hence in particular it takes trivial fibrations between fibrant objects to weak equivalences. Therefore, by Ken Brown's lemma, $U$ takes all weak equivalences between fibrant objects to weak equivalences. If $M\in \cM$ is a fibrant object, then $\eta^{\M}_M\colon M\to FM$ is a weak equivalence (a trivial cofibration, really) between fibrant objects and so $U(\eta^{\M}_M)=(U\eta^{\M})_M$ is a weak equivalence. In other words, the components of $U\eta^{\M}$ at fibrant objects are componentwise weak equivalences. The components of $CU\eta^{\M}\cdot\epsilon^{\cN} CU$ at fibrant objects are weak equivalences, by the 2-out-of-3 property. Hence $CU\eta^{\M}\cdot\epsilon^{\cN} CU$ is a componentwise homotopy equivalence when restricted to fibrant-cofibrant objects. Using the commutativity of the first diagram, we get  
\begin{align*}
CU\eta^{\M}\cdot\epsilon^{\cN} CU &= CUF\epsilon\cdot CUFV\epsilon^{\cN} U\cdot CU\eta^{\M} VCU\cdot C\eta CU \\
&= CU(F\epsilon\cdot FV\epsilon^{\cN} U)\cdot (CU\eta^{\M} V\cdot C\eta)CU,
\end{align*}
implying the latter is a componentwise homotopy equivalence at the level of fibrant-cofibrant objects. Taking $h_*$ of the last composition (always restricted to fibrant-cofibrant objects) we get that $\ou\overline{\epsilon}\cdot\overline{\eta}\ou$ is a componentwise equivalence. Finally, $\ou\epsilon'\cdot\eta' \ou=\ou\widetilde{\eta^{\M}}\cdot (\ou\overline{\epsilon}\cdot\overline{\eta}\ou)\cdot\widetilde{\epsilon^{\cN}}\ou$ by definition. Each piece of the composition is a componentwise equivalence, so the composition is also a componentwise equivalence. Hence, it is a pseudonatural equivalence, since a pseudonatural transformation is an equivalence if and only if it is a componentwise equivalence. The same argument applies to $\epsilon' \ov \cdot \ov\eta'$. Therefore Lemma \ref{lemma biadjunction} gives us the required biadjunction.\end{proof}

\begin{cor}\label{biadjoint biequivalence}
If the Quillen adjunction in Proposition \ref{enriched Quillen induces biadj} is an equivalence, then the induced biadjunction is a biequivalence.
\end{cor}
\begin{proof}
Since the Quillen adjunction is an equivalence, the components of the unit and the counit at fibrant and cofibrant objects are equivalences and hence are sent to pseudonatural equivalences serving as units and counits of the induced biadjoint biequivalences.
\end{proof}

\subsection{Double pseudofunctoriality of the homotopy 2-category}\label{sec: double pseud}
We end the background section with a result due to Shulman (see \cite{shulman2011comparing}) that will be useful in \Cref{chap:higher der}. This is stated in the language of double categories. 
The relevant notion of functor in our context is that of a double functor which is pseudo in both
directions, which in the literature is called \emph{double pseudofunctor} (see \cite[\S 6]{shulman2011comparing}).
\begin{rmk}
Many concepts defined within a 2-category also make sense in the context of a double category. In particular, we have a notion of \emph{conjunction} (a double-categorical version of adjunction) and a notion of \emph{companions} (similar to 2-isomorphisms). Both companions and conjunctions give rise to a notion of mates. As it is expected, double (pseudo)functors preserve both of them (see \cite{shulman2011comparing}).
\end{rmk}
\begin{defn}
Let \underline{Model} be the double category with:
\begin{enumerate}
    \item \emph{objects}: model categories $\M$,
    \item \emph{horizontal arrows}: left Quillen functors $F\colon\M\xslashedrightarrow{}\M'$,
    \item \emph{vertical arrows}: right Quillen functors $U\colon\M\to\cN$,
    \item \emph{squares}: natural transformations of the form
    \[\begin{tikzcd}
	\M & {\M'} \\
	\cN & {\cN'}
	\arrow["F", "\shortmid"{marking}, from=1-1, to=1-2]
	\arrow["U"', from=1-1, to=2-1]
	\arrow["G"', "\shortmid"{marking}, from=2-1, to=2-2]
	\arrow["{U'}", from=1-2, to=2-2]
	\arrow["\alpha", shorten <=12pt, shorten >=12pt, Rightarrow, from=2-1, to=1-2]
\end{tikzcd}\]
\end{enumerate}
\end{defn}

The following proposition was proved by Shulman in \cite{shulman2011comparing}, and, as we will see later, it generalizes to the enriched setting.
\begin{prop}\label{dbl psfunctor h}
There exists a double pseudofunctor $\emph{h}\colon\underline{\emph{Model}} \to {\emph{Sq}(\Cat)}$ sending a model category to its homotopy category.
\end{prop}

We will now introduce an analogue of the double pseudofunctor discussed in Proposition \ref{dbl psfunctor h} for enriched model categories. In this context, we speak of enriched Quillen adjunctions, which are enriched adjunctions (see \cite[1.11]{kelly1982basic}) whose underlying adjunctions are Quillen ones.
\begin{defn}
Let $\sSet_{\joy}$-\underline{Model} be the double category with:
\begin{enumerate}
    \item \emph{objects}: $\sSet_{\joy}$-enriched model categories $\M$,
    \item \emph{horizontal arrows}: $\sSet_{\joy}$-enriched left Quillen functors $F\colon\M\xslashedrightarrow{}\M'$,
    \item \emph{vertical arrows}: $\sSet_{\joy}$-enriched right Quillen functors $U\colon\M\to\cN$,
    \item \emph{squares}: $\sSet_{\joy}$-enriched natural transformations of the form
    \[\begin{tikzcd}
	\M & {\M'} \\
	\cN & {\cN'}
	\arrow["F", "\shortmid"{marking}, from=1-1, to=1-2]
	\arrow["U"', from=1-1, to=2-1]
	\arrow["G"', "\shortmid"{marking}, from=2-1, to=2-2]
	\arrow["{U'}", from=1-2, to=2-2]
	\arrow["\alpha", shorten <=12pt, shorten >=12pt, Rightarrow, from=2-1, to=1-2]
\end{tikzcd}\]
\end{enumerate}
\end{defn}

\begin{defn}
Let $\h_2\Gray$ be the $2$-category with:
\begin{enumerate}
\item \emph{objects}: $2$-categories $\C$,
\item \emph{1-cells}: strict $2$-functors $G\colon\C\to\D$,
\item \emph{2-cells}: equivalence classes of pseudonatural transformations $[\alpha]\colon G\Rightarrow G'$, where $[\alpha]=[\beta]$ if and only if there exists an invertible modification $\Phi\colon\alpha\Rrightarrow\beta$ in $\Gray$ (seen as self-enriched).
\end{enumerate}
\end{defn}
The following proposition will be useful in the next section to show that the example of a 2-derivator we are most interested in satisfies a version of (Der4), using the preservation of mates by double pseudofunctors.
\newpage
\begin{prop}\label{hocat double pseudofunctor enriched}
There exists a double pseudofunctor $$\widetilde{h_*}\colon\sSet_{\emph{Joyal}}\mbox{-}\underline{\emph{Model}}  \to{\emph{\Sq}(\h_2\Gray)}$$ such that
$$    \widetilde{h_*}(\M) = \h_*\M_{cf} 
$$
for every $\sSet_{\emph{Joyal}}$-enriched model category $\M$, with $\h_*$ being the change of base 2-functor induced by the homotopy category functor $h\colon\sSet\to\Cat$. 
\end{prop}
\begin{proof}
Define $\widetilde{h_*}$ on morphisms as follows:
\begin{align*}
    \widetilde{h_*}(G) &= \h_*F\circ\h_*G \\
    \widetilde{h_*}(U) &= \h_*C\circ\h_*U 
\end{align*}
for every $\sSet_{\joy}$-enriched left Quillen functor $G\colon\M\to\M'$ and every  $\sSet_{\joy}$-enriched right Quillen functor $U\colon\M\to\cN$, where $F$ and $C$ are the fibrant and cofibrant replacement simplicial functors of Definition \ref{simplicial replacements}. The action of $\widetilde{h_*}$ on a square 
\[\begin{tikzcd}
	\M & {\M'} \\
	\cN & {\cN'}
	\arrow["G", "\shortmid"{marking}, from=1-1, to=1-2]
	\arrow["U"', from=1-1, to=2-1]
	\arrow["{U'}", from=1-2, to=2-2]
	\arrow["\alpha", shorten <=12pt, shorten >=12pt, Rightarrow, from=2-1, to=1-2]
	\arrow["{G'}"', "\shortmid"{marking}, from=2-1, to=2-2]
\end{tikzcd}\]
sends it to 
\[\begin{tikzcd}
	{\widetilde{h_*}\M} & {\widetilde{h_*}\M'} \\
	{\widetilde{h_*}\cN} & {\widetilde{h_*}\cN'}
	\arrow["{\widetilde{h_*}G}", "\shortmid"{marking}, from=1-1, to=1-2]
	\arrow["{\widetilde{h_*}U}"', from=1-1, to=2-1]
	\arrow["{\widetilde{h_*}U'}", from=1-2, to=2-2]
	\arrow["{\widetilde{h_*}\alpha}", shorten <=12pt, shorten >=12pt, Rightarrow, from=2-1, to=1-2]
	\arrow["{\widetilde{h_*}G'}"', "\shortmid"{marking}, from=2-1, to=2-2]
\end{tikzcd}\]

\noindent in which $\widetilde{h_*}\alpha$ is the image through $h_*$ of the simplicial natural transformation whose components are the top rows of the lifting diagrams  
\[\begin{tikzcd}
	{(G'(UA)_c)_f} & {(U'(GA)_f)_c} \\
	{G'((UA)_c)} & {U'((GA)_f)} \\
	{G'(UA)} & {U'(GA)}
	\arrow["{\alpha_A}"', from=3-1, to=3-2]
	\arrow["U'i_{FA}"', from=3-2, to=2-2]
	\arrow["G'r_{UA}"', from=2-1, to=3-1]
	\arrow[from=2-1, to=2-2]
	\arrow["i_{G'(UA)_c}", "\sim"', hook, from=2-1, to=1-1]
	\arrow["r_{U'(GA)_f}", "\sim"', two heads, from=1-2, to=2-2]
	\arrow["\overline{\alpha}_A", dashed, from=2-1, to=1-2]
	\arrow["{\widetilde{\alpha}_A}", dashed, from=1-1, to=1-2]
\end{tikzcd}\]

The isomorphisms expressing the horizontal and vertical pseudofunctoriality of $\widetilde{h_*}$ are obtained in the same way as for $h$ of Proposition \ref{dbl psfunctor h}. Moreover the coherences and double naturality for squares give rise to the same diagrams of the aforementioned proposition. Therefore, everything we need to check for $\widetilde{h_*}$ holds by the exact same formal arguments that establish the double pseudofunctoriality of $\h\colon\underline{\text{Model}} \to {\text{Sq}(\Cat)}$ in Shulman's proof.
\end{proof}
\section{Higher derivator theory}\label{chap:higher der}
\graphicspath{{ch_other1/figures/}} 
\subsection{2-prederivators}
For a complete and cocomplete symmetric monoidal category $\V$, the functor category $[\A,\B]$ between $\V$-categories $\A$ and $\B$ is again $\V$-enriched
with homs given by the enriched coend
$$[\A,\B](T,S)=\int_{A\in\A}\B(TA,SA).$$
This is true in particular for $\V=\sSet$, so that $\sSet$-Cat is an object of $(\sSet$-Cat)-CAT. Recall from Example \ref{homotopy 2-cat} that the homotopy category functor 
$h\colon\sSet\to\Cat$ preserves finite products. Thus, we can define a functor 
$h_*\colon\sSet$-Cat$\to\twoCat$, which also preserves finite products, giving rise in turn to a 2-functor
\begin{center}
    $h_{**}\colon(\sSet$-Cat)-CAT$\to(\twoCat)$-CAT.
\end{center} 
For this reason, we set 
\begin{defn}
    $\Dia\coloneqq h_{**}(\sSet$-Cat).
\end{defn}
This is the $3$-category with objects small simplicially enriched categories and homs $\Dia(\A,\B)$ defined to be $h_*[\A,\B]$, the homotopy $2$-category of the simplicially enriched category $[\A,\B]$. Spelling this out:
\begin{itemize}
    \item[i)] the 1-cells of $\Dia$ are simplicially enriched functors,
    \item[ii)] the 2-cells of $\Dia$ are simplicially enriched natural transformations,
    \item[iii)] the 3-cells of $\Dia$ are homotopy classes of formal composites of modifications (i.e.\ 2-cells in $[\A,\B]$).
\end{itemize}
 The fully faithfulness of the nerve functor allows us to recover the previous definition of the higher category of diagrams given in \cite{NDVmres}, since the full sub-3-category of $\Dia$ on small 2-categories is the usual 3-category of 2-categories.

Before defining 2-prederivators, let us briefly review what Gray-categories are, along with some basic facts about them.

\begin{defn}
A $\Gray$-category is a category enriched over the monoidal category $\Gray=(\twoCat_0, \otimes, \bb1)$ of $2$-categories and strict $2$-functors equipped with the pseudo Gray tensor product (see e.g.\ \cite{JohYau:2-cat, street1988gray} for an explicit construction of the Gray tensor product of 2-categories).
\end{defn}
\begin{rmk}
    The monoidal category $\Gray$ is both symmetric and closed \cite{bourke2017gray}. If $\A,\B\in\Gray$, we will denote their internal hom with respect to this structure as $[\A,\B]_p$. This is the 2-category of strict 2-functors, \emph{pseudonatural transformations} and modifications between $\A$ and $\B$. The self-enrichment of $\twoCat_0$ over the pseudo Gray tensor product with internal homs given by $[-,-]_p$ is the $\Gray$-category $\Gray$.
\end{rmk}
\begin{defn}\label{quotient map from Gray tp to cartesian product}
    For a pair $\A,\B\in\Gray$, $[\A,\B]$ is the wide, locally full sub-2-category of strict 2-natural transformations in $[\A,\B]_p$. By the adjunction between the cartesian (resp.\ Gray tensor) product and the strict (resp.\ pseudo) hom, the inclusion $[\A,\B]\hookrightarrow[\A,\B]_p$ induces a monoidal functor 
    $$(\twoCat_0, \times, \bb1)\to (\twoCat_0, \otimes, \bb1),$$
 whose underlying functor is the identity functor on $\twoCat_0$. This monoidal functor preserves the unit since both monoidal structures have the same unit. Given a pair of 2-categories $\A,\B$, the associated comparison 2-functor
 $\A\otimes\B\to\A\times\B$ maps each generating interchange 
\[\begin{tikzcd}
	{A\otimes B} && {A'\otimes B} \\
	{A\otimes B'} && { A'\otimes B'}
	\arrow["{f\otimes \id_B}", from=1-1, to=1-3]
	\arrow["{\id_A\otimes g}"', from=1-1, to=2-1]
	\arrow["{\id_{A'}\otimes g}", from=1-3, to=2-3]
	\arrow["\cong"{description}, "f\otimes g", draw=none, from=2-1, to=1-3]
	\arrow["{f\otimes \id_{B'}}"', from=2-1, to=2-3]
\end{tikzcd}\]
in $\A\otimes\B$ to the identity square 
\[\begin{tikzcd}
	{(A, B)} && {(A',B)} \\
	{(A,B')} && { (A',B')}
	\arrow["{(f, \id_B)}", from=1-1, to=1-3]
	\arrow["{(\id_A, g)}"', from=1-1, to=2-1]
	\arrow["{(\id_{A'}, g)}", from=1-3, to=2-3]
	\arrow["{(f,\id_{B'})}"', from=2-1, to=2-3]
\end{tikzcd}\]
 in $\A\times\B$. Coherences for the unit are immediate, while associativity requires more work but it's fairly routine so we omit it. This implies that every 3-category may canonically be regarded as being a $\Gray$-category. A $\Gray$-category arises from a 3-category in this way if and only if each of its middle 4 interchanges is an identity 3-cell.
\end{defn}
\begin{defn}\label{2-prederivator}
A \emph{2-prederivator} is a \emph{trihomomorphism} $\bbD\colon\Dia^{\op}\to\GRAY$. 
\end{defn}
For a general definition of trihomomorphism between tricategories we refer to \cite[Definition 3.3.1]{gurski2006algebraic}. In analogy with ordinary derivator theory, we will use the superscript $(- )^*$ to denote the image through a 2-prederivator $\bbD$ of a $n$-cell for $n=1,2,3$.

\begin{exmp}
For every $\C\in\Dia$, we can consider the \emph{represented 2-prederivator} $\Dia(-,\C)\colon\Dia^{\op}\to\GRAY$, for which every coherence is an identity.
\end{exmp}    
\begin{exmp}
    The category of simplicially enriched categories is monoidal with respect to the product inherited by the categorical product at the level of hom-simplicial sets, and the homotopy 2-category 2-functor preserves products. Therefore, for every $\A\in\Dia$ and every 2-prederivator $\bbD$ we can define the \emph{shifted 2-prederivator} $\bbD^{\A}\coloneqq\bbD\circ(-\times\A).$
\end{exmp}
The key example of higher (pre)derivator is the one associated to a combinatorial\footnote{That is, a $\sSet_{\joy}$-enriched model category whose underlying model category is combinatorial.} $\sSet_{\joy}$-model category $\M$. For every small $2$-category $\J$, we considered the diagram category $[\J,\M]$. The latter is again $\sSet_{\joy}$-enriched and admits both the projective and the injective enriched model structure, which are equivalent by Proposition \ref{inj equiv proj}. This was done using a different and weaker notion of 2-prederivator. In the following we will show that the assignment
\begin{align*}   
    \bbD_{\M}\colon\Dia^{\op}&\longrightarrow\GRAY \\
                       \I&\mapsto h_*[\I,\M]_{cf}^{\proj}\\
                       (\I\xrightarrow{g}\J)&\mapsto h_*([\J,\M]_{cf}^{\proj}\xrightarrow{-\circ g}[\I,\M]_f^{\proj}\xrightarrow{C}[\I,\M]_{cf}^{\proj})
                       \end{align*}
defines a 2-prederivator for this new choice of $\Dia$, in which $\I$ is allowed to be a more generic simplicially enriched category. But first, let us take a brief detour to examine how the shifts of the 2-prederivators we have introduced so far look.
\begin{rmk}
    The represented 2-prederivator on $\C\in\Dia$ shifted by $\A\in\Dia$ is the 2-prederivator represented by $[\A, \C]$. This comes directly from the product-internal hom adjunction for simplicial categories. Indeed we have
    \begin{align*}
        \Dia(-,\C)^{\A} &= \Dia(-\times\A,\C)\\
        &\cong\Dia(-,[\A,\C]).
    \end{align*}
\end{rmk}
\begin{rmk}\label{shift of 2-pder of an enriched model category}
    Given small 2-categories $\A$, $\C\in\Dia$ and a  combinatorial $\sSet_{\joy}$-model category $\M$, there is an enriched Quillen equivalence (which is an isomorphism between the enriched categories that carry the model structures) $$[\A\times\C,\M]^{\proj}\simeq [\C,[\A,\M]^{\proj}]^{\proj}$$ coming from the product-internal hom adjunction plus the fact that weak equivalences are defined componentwise. Therefore    
    $$[\A\times\C,\M]_{cf}^{\proj}\simeq [\C,[\A,\M]^{\proj}]_{cf}^{\proj}$$
    that $h_*$ sends to an isomorphism
    $$h_*[\A\times\C,\M]_{cf}^{\proj}\cong h_*[\C,[\A,\M]_{cf}^{\proj}]^{\proj}$$ hence the 2-prederivator\footnote{By this we mean the shift of $\bbD_{\M}$ by $\A$. At this stage, we have yet to prove that $\bbD_{\M}$ is a 2-prederivator, so this definition is not yet meaningful—but it will be shortly.} $\bbD_{\M}^{\A}\coloneqq (\bbD_{\M})^{\A}$ is equivalent to $\bbD_{[\A,\M]^{\proj}}$. 
\end{rmk}
To prove that $\bbD_{\M}$ is actually a 2-prederivator, i.e.\ it is a trihomomorphism, we will proceed in a series of steps. First of all we will prove the following lemma.
\begin{lemma}\label{3-functor part of D_M}
For any simplicially enriched category $\A$ and any combinatorial $\sSet_{\joy}$-model category $\M$, the assignment
$$ \A\mapsto h_*[\A,\M]^{\proj}_f$$ defines a $3$-functor $\Dia^{\op}\to\twoCat$.
\end{lemma}
\begin{proof}
    Let us define $\qModr\in(\sSet\mbox{-}\Cat)\mbox{-}\CAT$ having combinatorial $\sSet_{\joy}$-model categories as objects and such that $\qModr(\M,\cN)$ is the full sub-simplicially enriched category of the simplicially enriched category $\sSet\mbox{-}\Cat(\M,\cN)$ spanning the enriched right Quillen functors $U\colon\M\to\cN$. Horizontal composites are defined as in $\sSet\mbox{-}\Cat(\M,\cN)$. We can then obtain $[-,\M]_f^{\proj}$ as the composite
$$\sSet\mbox{-}\Cat^{\op}\xrightarrow{[-,\M]^{\proj}}\qModr\xrightarrow{(-)_{f}}\sSet\mbox{-}\CAT$$
where the first component is obtained from the $(\sSet\mbox{-}\Cat)$-functor $$[-,\M]\colon\sSet\mbox{-}\Cat^{\op}\to\sSet\mbox{-}\CAT$$ by restricting the codomain to $\qModr$. We define $(-)_f$ to be the enriched functor that takes as input a combinatorial $\sSet_{\joy}$-model category and restricts it to the full subcategory of fibrant objects.
    
We move now to $h_{*}\colon\sSet\mbox{-}\CAT\to\twoCat.$ To find it, let us consider the functor
$$h_*\colon\sSet\mbox{-}\Cat_0\to\twoCat_0$$
between the underlying categories. We want to enrich it to a $(\sSet\mbox{-}\Cat)$-functor. To do this, we have to define simplicially enriched functors
$$[\A,\B]\to[h_*\A,h_*\B]$$ for every pair of simplicially enriched categories $\A$ and $\B$. Observe that if $\A$ and $\B$ are simplicially enriched categories,  $h_*\A$ and $h_*\B$ are 2-categories so that $[h_*\A,h_*\B]$ is again a 2-category (because $\twoCat$ is a 3-category), therefore $$[h_*\A,h_*\B]\cong N_{*}[h_*\A,h_*\B]$$ where we see $[h_*\A,h_*\B]$ on the RHS as a 2-category. Observe that we have a chain of natural bijections
\begin{prooftree}
\AxiomC{$[\A,\B]\to N_{*}[h_*\A,h_*\B]$}
\UnaryInfC{$h_*[\A,\B] \to [h_*\A,h_*\B]$}
\UnaryInfC{$h_*[\A,\B]\times h_*\A\to h_*\B$}
\end{prooftree}
where we used $h_*\dashv N_*$ and the product-internal hom adjunction. Defining the bottom morphism as $h_{*}(\ev_{\B})\colon h_*([\A,\B]\times\A)\cong h_*[\A,\B]\times h_*\A\to h_*\B$ and moving up the chain of equivalences we get the simplicially enriched functors we were looking for.

Considering the composite
$$\sSet\mbox{-}\Cat^{\op}\xrightarrow{[-,\M]_f^{\proj}}\sSet\mbox{-}\CAT\xrightarrow{h_*}\twoCat$$
in ($\sSet$-Cat)-CAT, we then conclude the existence of the desired 3-functor using the universal property of the unit of the adjunction $h_{**}\dashv N_{**}$
\[\begin{tikzcd}
	{\sSet\mbox{-}\Cat^{\op}} && {N_{**}\Dia^{\op}=N_{**}h_{**}(\sSet\mbox{-}\Cat^{\op})} \\
	\\
	&& {N_{**}(\twoCat)}
	\arrow["{\eta_{\sSet\mbox{-}\Cat^{\op}}}", from=1-1, to=1-3]
	\arrow[from=1-1, to=3-3]
	\arrow[dashed, from=1-3, to=3-3]
\end{tikzcd}\]
and the fully faithfulness of $N_{**}$.    
\end{proof}
\begin{rmk}
    If $\M$ is a combinatorial $\sSet_{\joy}$-model category, then $h_*(\M_{cf})$ is a \emph{bi-coreflective} full sub-2-category of $h_*\M_f$ in $\GRAY$. Specifically, in Definition \ref{simplicial replacements} we introduced a $\sSet_{\joy}$-enriched cofibrant replacement functor $C^{\M}$ and an enriched natural trasformation $\epsilon^{\M}$, whose images through $h_*$ give us the following data (we omit $h_*$ in the image for better readability):
\begin{itemize}
    \item a bi-coreflection 2-functor, i.e.\ a right biadjoint 
\[\begin{tikzcd}
	{h_*\M_f} && {h_*\M_{cf}}
	\arrow[""{name=0, anchor=center, inner sep=0}, "{C^{\M}}", curve={height=-18pt}, from=1-1, to=1-3]
	\arrow[""{name=1, anchor=center, inner sep=0}, "I^{\M}", curve={height=-18pt}, hook', from=1-3, to=1-1]
	\arrow["{\dashv^{b}}"{marking, allow upside down}, draw=none, from=1, to=0]
\end{tikzcd}\]
where $C^{\M}(\A)$ is a cofibrant replacement $\A_c$ of $\A$ in the model category $\M$.
\item A counit 
$$\epsilon^{\M}\colon I^{\M}C^{\M}\Rightarrow\id_{\M_f}$$
which is 2-natural and whose components $\epsilon_{\A}^{\M}\colon\A_c\to\A$ are trivial fibrations in $\M$. 
\item A unit 
$$\eta^{\M}\colon\id_{\M_cf}\Rightarrow C^{\M}I^{\M}$$
which is pseudonatural. This is constructed by a lifting argument which ensures that 
$$(\epsilon^{\M}I^{\M})\cdot(I^{\M}\eta^{\M}) = \id_{I^{\M}}.$$ It should also be noted that, since $\eta^{\M}$ is the unit of a bi-coreflection it is, in particular, a (left adjoint) equivalence.
\item A triangle isomorphism
$$\Phi^{\M}\colon (C^{\M}\epsilon^{\M})\cdot(\eta^{\M}C^{\M}) \cong \id_{\C^{\M}}$$
\end{itemize}
\end{rmk}

\begin{rmk}
    It is a basic result of $\Gray$-category theory, originally proved in Verity's thesis \cite{verity2011enriched}, that we may select the triangle isomorphism $\Phi^{\M}$ to satisfy two extra \emph{swallowtail identities}, one of which is the following:
\[\begin{tikzcd}
	& IC \\
	& ICIC \\
	IC && IC && {\id_{\epsilon}} \\
	& {\id_{\M_f}}
	\arrow["{I\eta C}", from=1-2, to=2-2]
	\arrow[""{name=0, anchor=center, inner sep=0}, curve={height=12pt}, Rightarrow, no head, from=1-2, to=3-1]
	\arrow[""{name=1, anchor=center, inner sep=0}, curve={height=-12pt}, Rightarrow, no head, from=1-2, to=3-3]
	\arrow["{\epsilon I C}", from=2-2, to=3-1]
	\arrow["{IC\epsilon}"', from=2-2, to=3-3]
	\arrow["{=}"{description}, draw=none, from=3-1, to=3-3]
	\arrow["\epsilon"', from=3-1, to=4-2]
	\arrow["{=}"{description}, draw=none, from=3-3, to=3-5]
	\arrow["\epsilon", from=3-3, to=4-2]
	\arrow[shorten <=6pt, shorten >=6pt, Rightarrow, no head, from=0, to=2-2]
	\arrow["\cong"{description}, shift left=2, "I\Phi", draw=none, from=1, to=2-2]
\end{tikzcd}\]
where we dropped the superscript $\M$ for conciseness. The interchange appearing in the bottom square is an equality in $\GRAY$ because $\epsilon$ is 2-natural. In other words, one simply posits that $\epsilon\cdot (I\Phi) = \id_{\epsilon}$. The second swallowtail identity follows from this and is elided for brevity.
\end{rmk}
Our goal here is to prove the following theorem.
\begin{thm}\label{trihomo part of D_M}
    There exists a trihomomorphism $$h_*(-)_{cf}\colon\qModr_3\to\GRAY,$$ where $\qModr_3$ denotes the 3-category $h_{**}(\qModr)$. This maps each combinatorial $\sSet_{\joy}$-model category $\M$ to the homotopy 2-category $h_*\M_{cf}$ of its fibrant and cofibrant objects. 
\end{thm}
For brevity, we shall write $H$ for this trihomomorphism under construction. Our method will be to do as much of this construction as possible in the sub-3-category $\twoCat$, wherein the trihomomorphism axioms simplify substantially, rather than in the more complicated context of $\GRAY$ itself. This approach relies crucially on the fact that we have constructed the data witnessing the bicoreflection of $h_*\M_{cf}$ in $h_*\M_f$, for each combinatorial $\sSet_{\joy}$-model category $\M$, so that $I^{\M}, C^{\M}$ and $\epsilon^{\M}$ are all inside $\twoCat$.

\begin{proof}[proof of \Cref{trihomo part of D_M}]
    
Our construction is as follows:
\begin{enumerate}
    \item for $\M\in\qModr_3$, let $H(\M)\coloneqq h_*\M_{cf}$
    \item for $\M,\cN\in\qModr_3$, define the action 
    $$\qModr_3(\M,\cN)\xrightarrow{H^{\M,\cN}}\twoCat(h_*\M_{cf}, h_*\cN_{cf})$$
    by $H^{\M,\cN}(-)\coloneqq C^{\cN}\circ h_*(-)_f \circ I^{\M}.$ This is a well defined 2-functor, because $h_*(-)_f$ is a 3-functor acting on the hom
$\qModr_3(\M,\cN)$ which is the pre/post composition by the 2-functors $I^{\M}$ and $C^{\cN}$ in $\twoCat$.
\item weak functoriality: for each $\M,\cN,\cP\in\qModr_3$, we define a 2-natural transformation 
\[\begin{tikzcd}
	{\qModr_3(\cN,\cP)\times\qModr_3(\M,\cN)} && {\qModr_3(\M,\cP)} \\
	\\
	{\twoCat(h_*\cN_{cf},h_*\cP_{cf})\times\twoCat(h_*\M_{cf},h_*\cN_{cf})} && {\twoCat(h_*\M_{cf},h_*\cP_{cf})}
	\arrow["\circ", from=1-1, to=1-3]
	\arrow["{H^{\cN,\cP}\times H^{\M,\cN}}"', from=1-1, to=3-1]
	\arrow["{H^{\M,\cP}}", from=1-3, to=3-3]
	\arrow["{\chi^{\M,\cN,\cP}}", shorten <=54pt, shorten >=54pt, Rightarrow, from=3-1, to=1-3]
	\arrow["\circ"', from=3-1, to=3-3]
\end{tikzcd}\]
Given a pair of right Quillen functors $U\colon\M\to\cN$ and $V\colon\cN\to\cP$ in $\qModr_3$, the component $\chi^{\M,\cN,\cP}_{U,V}$ is defined using the counit $\epsilon^{\cN}\colon I^{\cN}C^{\cN} \Rightarrow\id_{\cN_f}$, which by construction lies in $\twoCat$ (notice that $-$ with an abuse of notation $-$ we identify whiskering with a 1-cell and horizontal composition with the identity 2-cell on that 1-cell):
\[\begin{tikzcd}
	&&&& {H(VU)} \\
	{H(V)\circ H(U)} &&&& {C^{\cP}\circ h_{*}(VU)_f\circ I^{\M} } \\
	{\C^{\cP}\circ h_*V_f\circ I^{\cN}\circ C^{\cN}\circ h_*U_f \circ I^{\M}} &&&& {\C^{\cP}\circ h_*V_f\circ h_*U_f \circ I^{\cM}}
	\arrow[shorten <=3pt, shorten >=3pt, Rightarrow, no head, from=1-5, to=2-5]
	\arrow[shorten <=3pt, shorten >=3pt, Rightarrow, no head, from=2-1, to=3-1]
	\arrow[shorten <=3pt, shorten >=3pt, Rightarrow, no head, from=2-5, to=3-5]
	\arrow["{\C^{\cP}\circ h_*V_f\circ \epsilon^{\cN}\circ h_*U_f \circ I^{\M}}", from=3-1, to=3-5]
\end{tikzcd}\]
where in the codomain we used the 3-functoriality of $h_*(-)_f$. 

It is now a matter of simple calculations in the 3-category $\twoCat$ to show that 
\begin{itemize}
    \item[a)] the components $\chi_{U,V}\colon H(V) \circ H(U)\to H(VU)$ are 2-natural in $U$ and $V$,
    \item[b)] for $U\colon\M\to\cN, V\colon\cN\to\cP, W\colon\cP\to\cR$ the following square commutes
\[\begin{tikzcd}
	{H(W)\circ H(V)\circ H(U)} && {H(W)\circ H(VU)} \\
	{H(WV)\circ H(U)} && {H(WVU)}
	\arrow["{H(W)\circ\chi_{U,V}}", from=1-1, to=1-3]
	\arrow["{\chi_{W,V}\circ H(U)}"', from=1-1, to=2-1]
	\arrow["{\chi_{W,VU}}", from=1-3, to=2-3]
	\arrow["{=}"{description}, draw=none, from=2-1, to=1-3]
	\arrow["{\chi_{WV,U}}"', from=2-1, to=2-3]
\end{tikzcd}\]
\end{itemize}
\end{enumerate}

Setting aside, for the moment, the fact that we haven't shown that $\chi^{\M,\cN,\cP}$ is an equivalence, we can at least observe, using the last two points, that
\begin{enumerate}
    \item $\chi^{\M,\cN,\cP}$ is a 2-natural transformation which gives the composition comparison asked for in the third bullet point in Gurski's definition of trihomomorphism (see \cite[Definition 3.3.1.]{gurski2006algebraic}).
    \item b) tells us that the modification $\omega$ asked for in Gurski's fifth bullet point can be taken to be an identity. Notice also that the associativity squares labelled $a$ and $a'$ are identities in our case because we are working in a 3-category.
    \item Now to prove that the key 3-dimensional axiom relating instances of the modification $\omega$ $-$ which as we have shown in $2)$ may be taken to be an identity $-$ holds (i.e.\ Gurski's eighth bullet point) we may make the following observations:
    \begin{itemize}
        \item all modifications mentioned in Gurski's eighth bullet point that are, possibly whiskered, instances of $\omega$ are identities\footnote{This means that, in our setting, whenever we encounter the modification $\omega$ as defined by Gurski, or any modification obtained by whiskering it, it must be an identity due to the particular nature of our environment.}.
        \item All 2-cells labelled $a$ are associativity equivalences of the domain and codomain tricategories. In our case our tricategories are actually 3-categories, so these are actually identities.
        \item Similar comments apply to the modifications labelled $\pi$.
        \item The remaining, unlabelled, modifications are all also identities in our case. This is because they are structural isomorphisms of the domain/codomain tricategories, which in our case are identities since our tricategories are 3-categories. 
    \end{itemize}
So in our case, all modifications in Gurski's eighth bullet point are identities hence it holds trivially.
\end{enumerate}

Notice however that we haven't actually shown that $\chi^{\M,\cN,\cP}$ is an (adjoint) equivalence. Indeed it isn't an equivalence in $\twoCAT$ but it is an equivalence in $\GRAY$. To show this, fix $U\colon\M\to\cN, V\colon\cN\to\cP\in\qModr_3$ and consider the component
$$\chi_{U,V}\colon H(V)H(U)\to H(VU)$$
of $\chi^{\M,\cN,\cP}$. By definition, this is the 2-natural transformation $C^{\M}\circ h_*V_f\circ\epsilon^{\cN}\circ h_*U_f \circ I^{\M}.$ Now, given an object $A\in h_*\M_{cf}$ the component of this 2-natural transformation at $A$ is constructed by applying $C^{\M}\circ h_* V_f$ to the component of $\epsilon^{\cN}$ at $(h_*U_f\circ I^{\M})(A)=UA$, that being
$$C^{\M} V(\epsilon^{\cN}_{UA})\colon C^{\M} V C^{\cN} U(A)\to C^{\M} VU(A),$$
where we dropped the instances of the inclusion functor for brevity. By construction, $\epsilon^{\cN}_{UA}$ is a trivial fibration whose codomain $U(A)$ is fibrant (since $A$ is fibrant and $U$ is right Quillen) and whose domain $C^{\cN} U(A)$ is both fibrant and cofibrant. Now consider the following naturality square
\[\begin{tikzcd}
	{C^{\M}V C^{\cN} U(A)} && {C^{\M}VU(A)} \\
	{VC^{\cN}U(A)} && {VU(A)}
	\arrow["{C^{\M}V\epsilon^{\cN}_{U(A)}}", from=1-1, to=1-3]
	\arrow["{\epsilon^{\M}_{VC^{\cN}U(A)}}"', "\sim", two heads, from=1-1, to=2-1]
	\arrow["{\epsilon^{\M}_{VU(A)}}", "\sim"', two heads, from=1-3, to=2-3]
	\arrow["{V\epsilon^{\cN}_{U(A)}}"', from=2-1, to=2-3]
\end{tikzcd}\]
in which the upper horizontal arrow is the component of $\chi_{U,V}$ at $A$. Observe now that 
\begin{itemize}
    \item $U,V$ are right Quillen functors so they preserve fibrant objects, also $C^{\M}, C^{\cN}$ are cofibrant replacements so they take fibrant objects to ones that are fibrant and cofibrant. So all objects in the naturality square are fibrant and the upper two are also cofibrant.
    \item the components of $\epsilon^{\M}$ are all trivial fibrations, by construction, so the vertical morphisms in the naturality square are trivial fibrations as displayed.
    \item the component $\epsilon^{\cN}_{UA}$ is also a trivial fibration by construction and that property is preserved by the right Quillen $V$. Therefore the bottom horizontal arrow in the naturality square is a trivial fibration.
\end{itemize}
It follows by 2-out-of-3 for weak equivalences that the upper horizontal arrow of the naturality square is a weak equivalence between fibrant and cofibrant objects and hence that it is an equivalence in the 2-category $h_*\cP_{cf}$.
So we have shown that the 2-natural transformation
$$\chi_{U,V}\colon H(V)H(U)\to H(VU)$$
is a componentwise equivalence. Hence it has an equivalence inverse which is pseudonatural. That is, it is an equivalence in the 2-category $\GRAY(h_*\M_{cf}, h_*\cP_{cf})$. 

Thus, if we regard the 2-natural transformation $\chi^{\M,\cN,\cP}$ as a 2-cell with codomain 0-cell $\GRAY(h_*\M_{cf},h_*\cP_{cf})$ rather than its subcategory $\twoCAT(h_*\M_{cf}, h_*\cP_{cf})$, it is a componentwise equivalence and consequently is an equivalence in $\GRAY$. Of course, we may regard the structure $H$ that we have defined before as landing in $\GRAY$ rather than $\twoCAT$. By doing so, the coherence properties demonstrated so far remain unchanged but now we know that each $\chi^{\M,\cN,\cP}$ is a (left adjoint) equivalence as required of a trihomomorphism.

It is left to show that identities are weakly preserved. This means that, for each $\M\in\qModr_3$ we must also provide a (left adjoint) equivalence
\[\begin{tikzcd}
	\bb1 & {\qModr_3(\M,\M)} \\
	\\
	& {\GRAY(h_*\M_{cf},h_*\M_{cf})}
	\arrow["{\id_{\M}}", from=1-1, to=1-2]
	\arrow[""{name=0, anchor=center, inner sep=0}, "{\id_{h_*\M_{cf}}}"', from=1-1, to=3-2]
	\arrow["{H^{\M,\M}}", from=1-2, to=3-2]
	\arrow["{\iota^{\M}}"', shorten <=11pt, shorten >=8pt, Rightarrow, from=0, to=1-2]
\end{tikzcd}\]
as in Gurski's fourth bullet point. Applying the definition of $H^{\M,\M}$ above we see that the upper right composition of the triangle above maps the unique object $*\in\bb1$ to the 1-cell $C^{\M}\circ h_*(\id_{\M})_f \circ I^{\M}= C^{\M}\circ I^{\M},$ since $h_*(-)_f$ is a 3-functor so $h_*(\id_{\M})_f=\id_{h_*\M_f}$. It follows that to construct $\iota^{\M}$ all we need to do is to provide a pseudonatural (left adjoint) equivalence
$\id_{h_*\M_{cf}}\Rightarrow C^{\M}\circ I^{\M}$ and the obvious choice here is the unit $\eta^{\M}$. We know this is a (left adjoint) equivalence because it is the unit of a bi-coreflection. It is now easily checked that the modifications in Gurski's sixth bullet point may be taken to be identities. This follows directly by application of the triangle identity $(\eta^{\M}I^{\M})\cdot(I^{\M}\eta^{\M}) = \id_{I^{\M}}.$ To complete our construction it remains to show that Gurski's ninth bullet point holds. That follows again by arguing much as we did for Gurski's sixth bullet point to show that all modifications in those pastings are in fact identities in this case.
\end{proof}

 The previous results combined give a proof of the following theorem.
\begin{thm}
    If $\M$ is a combinatorial $\sSet_{\joy}$-enriched model category, $\bbD_{\M}$ is a 2-prederivator.
\end{thm}
\begin{proof}
For $\I\in\Dia$ and $\M$ a combinatorial $\sSet_{\joy}$-enriched model category, let us consider the composite assignment  
$$\I\mapsto h_*[\I,\M]_f^{\proj}\mapsto h_*[\I,\M]^{\proj}_{cf}.$$
We have seen from Lemma \ref{3-functor part of D_M} how to extend the first assignment to a 3-functor and from \Cref{trihomo part of D_M} how to extend the second one to a trihomomorphism. Composing these two, we get the trihomomorphism that we were looking for.
\end{proof}

\subsection{2-derivators}
In this section we introduce some of the axioms that a 2-prederivator should satisfy to deserve the name \emph{2-derivator}, inspired both by ordinary derivator theory and the theory of $\infty$-cosmoi. We propose to denote the axioms (HDer $n$) after their 1-derivatorial counterparts, where ``H'' stands for ``Higher''. 
\subsubsection{Axiom on (co)products}
The first axiom is a direct generalisation of (Der 1) from the theory of derivators. Notice that $\Dia$ is closed under coproducts.
\newline
\begin{mdframed}
\begin{axiom}[(HDer 1)]\label{axm: coprod}
   The canonical map $\bbD(\I\sqcup\J)\to\bbD(\I)\times\bbD(\J)$ is an equivalence for every $\I,\J\in\Dia$. In addition,  $\bbD(\varnothing)\simeq\bb1$.
\end{axiom}
\end{mdframed}
\begin{prop}
    The represented 2-prederivator satisfies \ref{axm: coprod}.
\end{prop}
\begin{proof}
    Being a representable (enriched) functor, the represented $2$-prederivator $\Dia(-,\C)$ sends weighted colimits to weighted limits. In particular, considering the coproduct and product as enriched limits\footnote{Note that this proof relies on the fact that these coproducts are colimits in $\Dia$, not merely in $\sSet\mbox{-}\Cat$. While not all colimits in $\sSet\mbox{-}\Cat$ are also colimits in $\Dia$, products and coproducts are.} we have that $\Dia(\I\sqcup\J, \C)\cong\Dia(\I,\C)\times\Dia(\J,\C)$ for every $\C$. 
\end{proof}
\begin{prop}
    If a 2-prederivator $\bbD$ satisfies \ref{axm: coprod}, then the shifted 2-prederivator $\bbD^{\A}$ satisfies \ref{axm: coprod} for every $\A\in\Dia$.
\end{prop}
\begin{proof}
    For every $\A\in\Dia$, the simplicially enriched functor $-\times\A$ is an enriched left adjoint by \cite[\S  2.3]{kelly1982basic}, therefore it preserves enriched coproducts.
    If $\bbD$ satisfies \ref{axm: coprod}, it follows that
\begin{align*}
\bbD^{\A}(\I\sqcup\J)&=\bbD((\I\sqcup\J)\times\A)\\
&\cong\bbD((\I\times\A)\sqcup(\J\times\A))\\
&\cong\bbD^{\A}(\I)\times\bbD^{\A}(\J). \qedhere
\end{align*}
\end{proof}
\begin{prop}
    The 2-prederivator $\bbD_{\M}$ associated to a combinatorial $\sSet_{\joy}$-enriched model category $\M$ satisfies \ref{axm: coprod}. 
\end{prop}
\begin{proof}
    Since $h$ commutes with finite products, it is enough to show that the isomorphism $$[\C\sqcup\D,\M]_{cf}^{\proj}\cong[\C, \M]_{cf}^{\proj}\times[\D,\M]_{cf}^{\proj}$$
    is also a Quillen equivalence between the corresponding enriched projective model structures. In order to prove the Quillen equivalence above, just notice that this isomorphism identifies fibrations and weak equivalences in the LHS with a pair of fibrations and weak equivalence in the RHS serving as components of the formers. On the other hand, cofibrations are completely determined by their left lifting property with respect to trivial fibrations. For instance, in $[\C, \M]_{cf}^{\proj}\times[\D,\M]_{cf}^{\proj}$ a morphism is a cofibration if and only if the following lifting problem
  \[\begin{tikzcd}
	{(H,K)} & {(F,G)} \\
	{(H',K')} & {(F',G')}
	\arrow[from=1-1, to=1-2]
	\arrow["f"', shift right=2, from=1-1, to=2-1]
	\arrow[from=2-1, to=2-2]
	\arrow["q", "\sim"'{sloped}, shift left=3, two heads, from=1-2, to=2-2]
	\arrow[dashed, from=2-1, to=1-2]
	\arrow["g", shift left=3, from=1-1, to=2-1]
	\arrow["p", "\sim"'{sloped}, shift right=2, two heads, from=1-2, to=2-2]
\end{tikzcd}\]
has a solution. We can then go back to $[\C\sqcup\D,\M]_{cf}^{\proj}$ assuming each time $p=\id_F$ or $q=\id_G$. Therefore $\bbD_{\M}(\C\sqcup\D)\cong\bbD_{\M}(\C)\times\bbD_{\M}(\D)$. \qedhere
\end{proof}
\subsubsection{Axiom on componentwise equivalences} 
In what follows, we provide a formalization of the idea that equivalences are detected componentwise. As we saw in \Cref{chap:background}, the underlying diagram functor for a prederivator is defined using the product-internal hom adjunction available in $\Cat$. In a similar way, one can define \emph{underlying diagram 2-functors} for a 2-prederivator.
\begin{defn}
    For any $\C\in\Dia$, the underlying diagram 2-functor $$\dia_{\C}\colon\bbD(\C)\to[h_*\C,\bbD(\bb1)]_p$$ is the output of the chain of transpositions 
    \begin{prooftree}
\AxiomC{$h_*\C\cong\Dia(\bb1,\C) \xrightarrow{\bbD_{\bb1,\C}}[\bbD(\C),\bbD(\bb1)]_p$}
\UnaryInfC{$\bbD(\C)\otimes h_*\C \to\bbD(\bb1)$}
\UnaryInfC{$\dia_{\C}\colon\bbD(\C)\to[h_*\C,\bbD(\bb1)]_p$}
\end{prooftree}
\end{defn}
\noindent where we used that $\Dia(\bb1,\C)=h_*[\bb1,\C]\cong h_*\C$, since $[\bb1,\C]\cong\C$. When $\C$ is a 2-category, the definition of the underlying diagram 2-functor simplifies to
$$\dia_{\C}\colon\bbD(\C)\to[\C,\bbD(\bb1)]_p$$
since $h_*\C\cong\C$ for every 2-category $\C$.
Remember that $[h_*\C,\bbD(\bb1)]_p$ is the internal hom in $\GRAY$ between $h_*\C$ and $\bbD(\bb1)$ so it is the 2-category of 2-functors $h_*\C\to\bbD(\bb1)$ (which, by adjunction, correspond to simplicially enriched functors $\C\to\bbD(\bb1)$), pseudonatural transformations between such 2-functors and modifications between those pseudonatural transformations.
\begin{rmk}
    The underlying diagram 2-functor acts on 0 and 1-cells like the underlying diagram functor from Definition \ref{underlying diagram functor}. The action on 2-cells is a direct generalisation of the one on 0- and 1-cells. 
\end{rmk}


\bigskip
\begin{mdframed}
    \begin{axiom}[(HDer 2)]\label{axm: cptwise equivalences}
     The 2-functors $$\emph{dia}_{\C}\colon\bbD(\C)\to[h_*\C,\bbD(\bb1)]_p$$ are conservative on $1$-cells for every $\C\in\Dia$.

\end{axiom}
\end{mdframed}

\noindent\textit{Notation.} We denote with $\dia_{\C}^{\A}\colon\bbD(\A\times\C)\to[h_*\C, \bbD(\A)]_p$ the underlying diagram 2-functor $\dia_{\C}$ relative to the shifted 2-prederivator $\bbD^{\A}$. In other words, it is the output of the chain of transpositions 

\begin{center}
    \begin{tabular}{rc}
    & \begin{tabular}{@{}c@{}}
\begin{tikzcd}
	{h_*\C\cong\Dia(\bb1,\C)} && {[\bbD(\A\times\C),\bbD(\A)]_p} \\
	& {\Dia(\A\times\bb1,\A\times\C)}
	\arrow[from=1-1, to=1-3]
	\arrow["{\A\times-}"'{pos=0.4}, from=1-1, to=2-2]
	\arrow["{\bbD_{\A,\A\times\C}}"', from=2-2, to=1-3]
\end{tikzcd}
    \end{tabular}
    \\
    \hline
    & $\bbD(\A\times\C)\otimes h_*\C\to\bbD(\A)$ \\
    \hline
    & $\bbD(\A\times\C) \xrightarrow{\dia_{\C}^{\A}}[h_*\C, \bbD(\A)]_p$
\end{tabular}
\end{center}

\bigskip

\begin{prop}
If a 2-prederivator $\bbD$ satisfies \ref{axm: cptwise equivalences}, then the shifted 2-prederivator $\bbD^{\A}$ satisfies \ref{axm: cptwise equivalences} for every $\A\in\Dia$.\end{prop}
\begin{proof}
    We want to prove that $\dia_{\C}^{\A}$ is conservative for every $\A$ and $\C$. Consider the following diagram 

\[\begin{tikzcd}
	{\bbD(\A\times\C)} & {[h_*\C,\bbD(\A)]_p} \\
	{[h_*\A\times h_*\C, \bbD(\bb1)]_p} & {[h_*\C, [h_*\A,\bbD(\bb1)]_p]_p}
	\arrow["{\dia_{\C}^{\A}}", from=1-1, to=1-2]
	\arrow["{\dia_{\A\times \C}}"', from=1-1, to=2-1]
	\arrow["\simeq"{description}, draw=none, from=1-1, to=2-2]
	\arrow["{[h_*\C, \dia_{\A}]_p}", from=1-2, to=2-2]
	\arrow[from=2-1, to=2-2]
\end{tikzcd}\]

\noindent commuting up to equivalence, where the map $$[h_*\A\times h_*\C, \bbD(\bb1)]_p\to[h_*\C, [h_*\A,\bbD(\bb1)]_p]_p\cong[h_*\A\otimes h_*\C, \bbD(\bb1)]_p$$ is induced by the quotient map $h_*\A\otimes h_*\C\to h_*\A\times h_*\C$ (see Definition \ref{quotient map from Gray tp to cartesian product}). The bottom map is conservative because the quotient map $h_*\A\otimes h_*\C\to h_*\A\times h_*\C$ is bijective on objects and the internal homs are with respect to the Gray tensor product, so equivalences are componentwise. Vertical maps are conservative by assumption (and in the case of the right map again because equivalences are componentwise). It follows that $\dia_{\C}^{\A}$ is conservative.
\end{proof}
\begin{prop}
    The 2-prederivator $\bbD_{\M}$ associated to a combinatorial $\sSet_{\joy}$-enriched model category $\M$ satisfies \ref{axm: cptwise equivalences}. 
\end{prop}
\begin{proof}
We need to show that the 2-functors 
    $$\dia_{\C}\colon h_*[\C,\M]^{\proj}_{cf}\to[h_*\C, h_*[\bb1,\M]^{\proj}_{cf}]_p \cong [h_*\C, h_*\M_{cf}]_p$$
     are conservative on 1-cells for every $\C\in\Dia$. This holds since the equivalences in the homotopy
2-category $h_{*}\M_{cf}$ are all weak equivalences in $\M$, and weak equivalences in the diagram categories are defined componentwise. 
\end{proof}
\newpage
\subsubsection{Axiom on biadjunctions}
The following axiom is a generalisation of axiom (Der 3). 
\newline
\begin{mdframed}
    \begin{axiom}[(HDer 3)]\label{axm: biadjunctions}
     Every $f\colon \A\to \B$ in $\Dia$ induces a biadjoint triple 
\[
\begin{tikzcd}
\bbD(\B)\ar[r,"f^*"description,""{name=A, below}] & \bbD(\A)\ar[l,bend right=49,"\Lan_f"',""{name=B,above}]\ar[l,bend left=49,"\Ran_f",""{name=C,above}] \ar[from=A, to=C, symbol=\perp_{b}, rotate=-90]\ar[from=B, to=A, symbol=\perp_{b}, rotate=-90]
\end{tikzcd}
\]
\end{axiom}
\end{mdframed}

We consider biadjunctions instead of 2-adjunctions, since our examples $-$ most notably 
$\bbD_{\M}$ $-$ are not strict in general. For instance, for every morphism $f$, the 2-functor $f^*$ is only defined up to equivalence. We will make this precise in the next result.
\begin{prop}\label{unique pullback}
Let $f\colon\D\to\C$ be in $\Dia$. The pullback $2$-functor $f^*\colon\bbD_{\M}(\C)\to\bbD_{\M}(\D)$ is uniquely determined up to pseudonatural equivalence. 
\end{prop}
\begin{proof}
We have to show that the following diagram (the vertical equivalences are a consequence of Corollaries \ref{equiv proj inj joyal enriched} and \ref{biadjoint biequivalence})
\begin{center}
    \begin{tikzcd}
h_*{{[\C,\M]}^{\proj}_{cf}} \ar[d, phantom, "\simeq"description, sloped, yshift=-2mm]\arrow[r, "-\circ f"]\arrow[d, "F^i"', bend right, xshift=-2mm] & h_*{{[\D,\M]}^{\proj}_{f}} \arrow[r, "C^p"] & h_*{{[\D,\M]}^{\proj}_{cf}} \ar[d, phantom, "\simeq", sloped, yshift=-2mm]\ar[dl, hook', bend right=8]\arrow[d, "F^i"', bend right, xshift=-2mm] \\
h_*{{[\C,\M]}^{\inj}_{cf}} \arrow[r, "-\circ f"'] \arrow[u, "C^p"', bend right, xshift=-2mm] & h_*{{[\D,\M]}^{\inj}_{c}} \arrow[r, "F^i"'] & h_*{{[\D,\M]}^{\inj}_{cf}} \arrow[u, "C^p"', bend right, xshift=-2mm] 
\end{tikzcd}
\end{center}
commutes up to a pseudonatural equivalence (we denote with $C^p,F^p$ and $C^i, F^i$ the cofibrant and fibrant replacements with respect to the projective and the injective model structure). The triangle on the right commutes strictly since every cofibration in the projective model structure is a cofibration in the injective model structure and hence every cofibrant object in the projective model structure is a cofibrant object in the injective model structure as well. A $2$-functor $g\colon\C\to\M$ in $h_*[\C,\M]^{\proj}_{cf}$ is sent to $(F^i g)f$ by $(-\circ f)F^i$, and it is sent to $C^p(gf)$ by the other composite arrow. We have a trivial fibration $C^p(gf)\xrightarrowdbl{\sim}gf$ in the projective model structure and a trivial cofibration $gf\xhookrightarrow{\sim}(F^ig)f$ in the injective model structure (since $g\xhookrightarrow{\sim}F^ig$ is a trivial cofibration and $-\circ f$ is left Quillen), so in particular they are levelwise weak equivalences. Hence the composite $C^p(gf) \xrightarrow{\sim}(F^ig)f$ is a pointwise weak equivalence between cofibrant objects. Then Corollary \ref{Cofib and fib replacement composed with pointwise equiv} implies that the image through $F^i$ of this morphism is the component in $g$ of a pseudonatural equivalence when we pass at the level of homotopy $2$-categories. Hence our diagram commutes up to pseudonatural equivalence meaning $f^*$ is uniquely determined up to pseudonatural equivalence.
\end{proof}

\begin{prop}
    The represented 2-prederivator $\Dia(-,\C)$ satisfies \ref{axm: biadjunctions} whenever $\C$ is complete and cocomplete as a $\sSet$-category. 
\end{prop}
\begin{proof}
This comes from enriched category theory: when $\C$ is complete and cocomplete we can construct left and right Kan extensions by means of weighted (co)limits \cite[\S 4.1]{kelly1982basic}. These adjoints lift to ($\sSet\mbox{-}\Cat$)-enrichment thanks to Kelly’s theorem on the lifting of adjunctions, applicable when the category admits certain colimits, and pass through $h_*$ for functoriality.
\end{proof}

\begin{prop}
 If a 2-prederivator $\bbD$ satisfies \ref{axm: biadjunctions}, then the shifted 2-prederivator $\bbD^{\A}$ satisfies \ref{axm: biadjunctions} for every $\A\in\Dia$.\end{prop}
    \begin{proof}
For every $f\colon\I\to\J$ in $\Dia$ there exists a biadjoint triple 
\[
\begin{tikzcd}
\bbD^{\A}(\J)\ar[r,"f^*"description,""{name=A, below}] & \bbD^{\A}(\I)\ar[l,bend right=49,"\Lan_f"',""{name=B,above}]\ar[l,bend left=49,"\Ran_f",""{name=C,above}] \ar[from=A, to=C, symbol=\perp_{b}, rotate=-90]\ar[from=B, to=A, symbol=\perp_{b}, rotate=-90],
\end{tikzcd}
\]
because the pullback 2-functor is given by the action of $\bbD$ on $\A\times f\colon\A\times\I\to\A\times\J$, which is still in $\Dia$, hence $\bbD^{\A}(f)=\bbD(\A\times f)$ admits homotopy left and right Kan extensions.
\end{proof}
\begin{prop}
    The 2-prederivator $\bbD_{\M}$ associated to a combinatorial $\sSet_{\joy}$-enriched model category $\M$ satisfies \ref{axm: biadjunctions}. 
\end{prop}
\begin{proof}
    Suppose we have a combinatorial $\sSet_{\joy}$-enriched model category $\M$ and $\C, \D\in\Dia$. We want to show that any $f\colon \D\to \C$  in $\Dia$ induces a biadjunction 
\[
\begin{tikzcd}
h_*[\C,\M]_{cf}\ar[r,"f^*"description,""{name=A, below}] & h_*[\D,\M]_{cf}\ar[l,bend right,"\text{lan}_f"',""{name=B,above}]\ar[l,bend left,"\text{ran}_f",""{name=C,above}] \ar[from=A, to=C, symbol=\perp_b, rotate=-85]\ar[from=B, to=A, symbol=\perp_b, rotate=-93].
\end{tikzcd}
\]
Using Proposition \ref{enriched Quillen induces biadj}, it is enough to show that the pullback functor $f^*$ sits inside an enriched adjoint triple, where the top adjoint pair is Quillen with respect to the projective model structure, and the bottom one is Quillen with respect to the injective model structure
\[
\begin{tikzcd}
{[\C,\M]}\ar[r,"f^*"description,""{name=A, below}] & {[\D,\M]}\ar[l,bend right,"\Lan_f"',""{name=B,above}]\ar[l,bend left,"\Ran_f",""{name=C,above}] \ar[from=A, to=C, symbol=\dashv, rotate=7]\ar[from=B, to=A, symbol=\dashv, rotate=-4].
\end{tikzcd}
\]
Since $\M$ is complete and cocomplete, we know from \cite[\S 4.3]{kelly1982basic} that the adjunctions $\Lan_f\dashv f^*$ and $f^*\dashv\Ran_f$ exist and are enriched. It remains to show that they are Quillen pairs, meaning that the adjunctions between the underlying categories are Quillen. From Theorem 5.4 in \cite{moser2019injective} we know that the underlying categories $[\C,\M]_0$ and $[\D,\M]_0$ admit both the projective and the injective model structure (which are again enriched) and furthermore Corollary \ref{equiv proj inj joyal enriched} ensures that these model structures are Quillen equivalent (so that we get equivalent homotopy $2$-categories), so let us choose them accordingly. We consider $$\Lan_f\dashv f^*\colon[\C,\M]_0^{\proj}\rightleftarrows[\D,\M]_0^{\proj} \ \text{and} \ f^*\dashv\Ran_f\colon[\C,\M]_0^{\inj}\rightleftarrows[\D,\M]_0^{\inj}$$ and prove that $f^*$ is both right and left Quillen with respect to the chosen model structures. Recall that $[\C,\M]_0$ and $[\D,\M]_0$ are categories with objects being simplicial functors $\C\to\M$ and $\D\to\M$ and morphisms between any pair of such functors $H,K$ being simplicial natural transformations $H\Rightarrow K$. Let us prove that $f^*$ is right Quillen, the other case being completely analogous. Notice that a simplicial natural transformation $\alpha\colon H\Rightarrow K$ in $[\C,\M]_0^{\proj}$ is a fibration (resp.\ a weak equivalence) if and only if $\alpha_c$ is a fibration (resp.\ a weak equivalence) in $\M_0$ for every $c\in\C$. Then $f^*\alpha\colon f^*H=Hf\Rightarrow f^*K=Kf$ has components $(f^*\alpha)_d = (\alpha f)_d = \alpha_{fd}$ for a generic $d\in\D$ and so it is a fibration or a weak equivalence whenever $\alpha$ is. In particular if $\alpha$ is a trivial fibration, then $f^*\alpha$ is a trivial fibration as well.
\end{proof}
\subsubsection{Axiom on exact squares}
All the equivalences appearing in squares like the one in the following definition are images of equalities of the kind $h = fg$ (i.e.\ compositions of 1-cells) that the 2-prederivator sends to equivalences - since it preserves composition of 1-cells up to equivalence. The actual name of those equivalences is not relevant for the proofs so we decided to omit them.
\begin{defn}
    A square in $\Dia$ 
\[\begin{tikzcd}
	A & B \\
	C & D
	\arrow["f", from=1-1, to=1-2]
	\arrow["k"', from=1-1, to=2-1]
	\arrow["g", from=1-2, to=2-2]
	\arrow["h"', from=2-1, to=2-2]
\end{tikzcd}\]
    is called \emph{$\bbD$-exact on the left} if its image through $\bbD$ 
\[\begin{tikzcd}
	{\bbD(D)} & {\bbD(C)} \\
	{\bbD(B)} & {\bbD(A)}
	\arrow["{h^*}", from=1-1, to=1-2]
	\arrow["{g^*}"', from=1-1, to=2-1]
	\arrow["\simeq"{description}, draw=none, from=1-1, to=2-2]
	\arrow["{k^*}", from=1-2, to=2-2]
	\arrow["{f^*}"', from=2-1, to=2-2]
\end{tikzcd}\]
    \noindent satisfies the \emph{left Beck-Chevalley} condition, i.e.\ the mate $$\Lan_f k^* \Rightarrow g^*\Lan_h$$ induced by the biadjunctions $\Lan_f\dashv_b f^*$ and $\Lan_h\dashv h^*$ is an equivalence. Similarly, we say that the square on the top is \emph{$\bbD$-exact on the right} if its image through $\bbD$ satisfies the \emph{right Beck-Chevalley} condition, i.e.\ the mate $$g^*\Ran_h \Rightarrow \Ran_f k^*$$ induced by the biadjunctions $f^*\dashv_b\Ran_f$ and $h^*\dashv_b\Ran_h$ is an equivalence. A square is \emph{$\bbD$-exact} if it is $\bbD$-exact on the left and $\bbD$-exact on the right. 
\end{defn}
\begin{rmk}
  Mates compose, meaning that given an isomorphism 
\[\begin{tikzcd}
	A & B & C & A & C \\
	X & Y & Z & X & Z
	\arrow[from=1-1, to=1-2]
	\arrow[""{name=0, anchor=center, inner sep=0}, from=1-1, to=2-1]
	\arrow[from=1-2, to=1-3]
	\arrow[""{name=1, anchor=center, inner sep=0}, from=1-2, to=2-2]
	\arrow[""{name=2, anchor=center, inner sep=0}, from=1-3, to=2-3]
	\arrow[from=1-4, to=1-5]
	\arrow[""{name=3, anchor=center, inner sep=0}, from=1-4, to=2-4]
	\arrow[""{name=4, anchor=center, inner sep=0}, from=1-5, to=2-5]
	\arrow[""{name=5, anchor=center, inner sep=0}, curve={height=-18pt}, from=2-1, to=1-1]
	\arrow["\alpha", shorten <=12pt, shorten >=12pt, Rightarrow, from=2-1, to=1-2]
	\arrow[from=2-1, to=2-2]
	\arrow[""{name=6, anchor=center, inner sep=0}, curve={height=-18pt}, from=2-2, to=1-2]
	\arrow["\beta", shorten <=12pt, shorten >=12pt, Rightarrow, from=2-2, to=1-3]
	\arrow[from=2-2, to=2-3]
	\arrow[""{name=7, anchor=center, inner sep=0}, curve={height=-18pt}, from=2-3, to=1-3]
	\arrow[""{name=8, anchor=center, inner sep=0}, curve={height=-18pt}, from=2-4, to=1-4]
	\arrow["\gamma", shorten <=12pt, shorten >=12pt, Rightarrow, from=2-4, to=1-5]
	\arrow[from=2-4, to=2-5]
	\arrow[""{name=9, anchor=center, inner sep=0}, curve={height=-18pt}, from=2-5, to=1-5]
	\arrow[symbol=\perp_b, rotate=94, draw=none, from=0, to=5]
	\arrow[symbol=\perp_b, rotate=94, draw=none, from=1, to=6]
	\arrow[symbol=\perp_b, rotate=94, draw=none, from=2, to=7]
	\arrow["\cong"{description}, draw=none, from=2, to=8]
	\arrow[symbol=\perp_b, rotate=94,draw=none, from=3, to=8]
	\arrow[symbol=\perp_b, rotate=94, draw=none, from=4, to=9]
\end{tikzcd}\]
in $\Gray$, we can paste the mate of the leftmost 2-cell along the specified biadjunctions with the mate of the second. In the middle of this pasting diagram we then have an instance of one of the triangular isomorphisms for the biadjunction involving $B$ and $Y$. Therefore we can cancel it out, and we get the mate of the 2-cell on the RHS up to isomorphism. 
\end{rmk}
Recall Definition \ref{collage}, where we defined the collage of a profunctor and that of a weight as a special case. Consider the weights $$\B(f-,b)\colon\A^{\op}\to\sSet \hspace{1cm}\text{and}\hspace{1cm}\B(b,f-)\colon\A\to\sSet$$ for $f\colon\A\to\B$ in $\Dia$ and $b\in\B$. 
\vspace{7mm}
\begin{mdframed}
    \begin{axiom}[(HDer 4)]\label{axm: ptwise Kan ext}
    For every $f\colon\A\to\B$ in $\Dia$ and every object $b\in\B$ the diagrams
\[\begin{tikzcd}
	\A & {\emph{coll}\B(f,b)} & \A & {\emph{coll}\B(b, f)} \\
	\A & \B & \A & \B
	\arrow["i_{\A}^{\B(f,b)}", hook, from=1-1, to=1-2]
	\arrow[Rightarrow, no head, from=1-1, to=2-1]
	\arrow["{\pi_{\B}}", from=1-2, to=2-2]
	\arrow["i_{\A}^{\B(b,f)}", hook, from=1-3, to=1-4]
	\arrow[Rightarrow, no head, from=1-3, to=2-3]
	\arrow["{\pi_{\B}}", from=1-4, to=2-4]
	\arrow["f"', from=2-1, to=2-2]
	\arrow["f"', from=2-3, to=2-4]
\end{tikzcd}\]
    are $\bbD$-exact. 
\end{axiom}
\end{mdframed}
\vspace{7mm}
First, we will show that if $\bbD$ satisfies \ref{axm: ptwise Kan ext}, then it satisfies a ``generalised elements'' version of it. We focus on the generalised version of the diagram on the left, since the other is analogous. 


\begin{prop}\label{generalised elements}
Suppose $\bbD$ satisfies \ref{axm: ptwise Kan ext} and consider a cospan $\A\xrightarrow{f}\C\xleftarrow{g}\B$ in $\Dia$. Then the diagram 
\[\begin{tikzcd}
	\A & {\emph{coll}(f,g)} \\
	\A & \C
	\arrow["{i_{\A}}", hook, from=1-1, to=1-2]
	\arrow[Rightarrow, no head, from=1-1, to=2-1]
	\arrow["{\pi_{\C}}", from=1-2, to=2-2]
	\arrow["f"', from=2-1, to=2-2]
\end{tikzcd}\]
is $\bbD$-exact.
\end{prop}
\begin{proof}
Notice that $\coll(i_{\A},i_{\B}b)=\coll(f,gb)$ since they have the same set of objects and 
\begin{align*}
    \coll(i_{\A},i_{\B}b)(a,\bullet)&=\coll(f,g)(i_{\A}a,i_{\B}b\bullet) \\
    &= \coll(f,g)(a,b) \\
    &=\C(fa,gb)\\
    &=\coll(f,gb)(a,\bullet).
\end{align*} 
Moreover there exists a simplicially enriched functor $j_{gb}\colon \coll(f,gb)\to\coll(f,g)$ corresponding to the inclusion of $gb$ in the collage. Considering the diagram
\[\begin{tikzcd}
	{\bbD(\coll(i_{\A},i_{\B} b))} & {\bbD(\coll(f,gb))} \\
	{\bbD(\A)} & {\bbD(\coll(f,g))} \\
	{\bbD(\A)} & {\bbD(\C)}
	\arrow[""{name=0, anchor=center, inner sep=0}, "{(i_{\A}^{gb})^*}"', from=1-1, to=2-1]
	\arrow["\simeq"', from=1-2, to=1-1]
	\arrow[Rightarrow, no head, from=2-1, to=3-1]
	\arrow[""{name=1, anchor=center, inner sep=0}, "{j_{gb}^*}"', from=2-2, to=1-2]
	\arrow[""{name=2, anchor=center, inner sep=0}, "{i_{\A}^*}", from=2-2, to=2-1]
	\arrow["{\pi_{\C}^*}"', from=3-2, to=2-2]
	\arrow[""{name=3, anchor=center, inner sep=0}, "{f^*}", from=3-2, to=3-1]
	\arrow["\simeq"{description}, draw=none, from=0, to=1]
	\arrow["\simeq", draw=none, from=2, to=3]
\end{tikzcd}\]
we see that the subdiagrams
\[\begin{tikzcd}
	{\bbD(\A)} && {\bbD(\coll(f,gb))} & {\bbD(\A)} && {\bbD(\coll(f,g b))} \\
	{\bbD(\A)} && {\bbD(\C)} & {\bbD(\A)} && {\bbD(\coll(f,g))}
	\arrow[""{name=0, anchor=center, inner sep=0}, "{\Lan_{i_{\A}^{gb}}}", curve={height=-24pt}, from=1-1, to=1-3]
	\arrow[Rightarrow, no head, from=1-1, to=2-1]
	\arrow[""{name=1, anchor=center, inner sep=0}, "{(i_{\A}^{gb})^*}", from=1-3, to=1-1]
	\arrow[""{name=2, anchor=center, inner sep=0}, "{\Lan_{i_{\A}^{gb}}}", curve={height=-24pt}, from=1-4, to=1-6]
	\arrow[Rightarrow, no head, from=1-4, to=2-4]
	\arrow["\simeq"{description}, draw=none, from=1-4, to=2-6]
	\arrow[""{name=3, anchor=center, inner sep=0}, "{(i_{\A}^{gb})^*}", from=1-6, to=1-4]
	\arrow["\simeq"{description}, draw=none, from=2-1, to=1-3]
	\arrow[""{name=4, anchor=center, inner sep=0}, "{\Lan_f}", from=2-1, to=2-3]
	\arrow["{j^*_{gb}\pi_{\C}^*}"', from=2-3, to=1-3]
	\arrow[""{name=5, anchor=center, inner sep=0}, "{f^*}", curve={height=-24pt}, from=2-3, to=2-1]
	\arrow[""{name=6, anchor=center, inner sep=0}, "{\Lan_{i_{\A}}}", from=2-4, to=2-6]
	\arrow["j_{gb}^*"', from=2-6, to=1-6]
	\arrow[""{name=7, anchor=center, inner sep=0}, "{i_{\A}^*}", curve={height=-24pt}, from=2-6, to=2-4]
	\arrow[symbol=\perp_b, rotate=65, yshift=1mm, draw=none, from=1, to=0]
	\arrow[symbol=\perp_b, rotate=-115, yshift=1mm, draw=none, from=2, to=3]
	\arrow[symbol=\perp_b, rotate=-90, yshift=1mm, draw=none, from=4, to=5]
	\arrow[symbol=\perp_b, rotate=-70, yshift=1mm, draw=none, from=6, to=7]
\end{tikzcd}\]

\noindent are both in a form for which we can apply \ref{axm: ptwise Kan ext}. We conclude using that mates compose and the 2-out-of-3 property of equivalences.
\end{proof}
\begin{prop}\label{Kan ext of fully faithful}
Suppose $\bbD$ satisfies \ref{axm: ptwise Kan ext} and let $I\colon\A\hookrightarrow\B$ be a fully faithful simplicial functor, meaning it is an isomorphism at the level of hom simplicial sets, then the unit $\eta\colon\id_{\bbD(\A)}\Rightarrow I^*\circ\Lan_I$ of the biadjunction
\[
\begin{tikzcd}
\bbD(\A)\ar[r,bend left,"\Lan_I",""{name=A, below}] & \bbD(\B)\ar[l,bend left,"I^*",""{name=B,above}] \ar[from=A, to=B, symbol=\perp_b, rotate=-90]
\end{tikzcd}
\]
is a pseudonatural equivalence. 
\end{prop}
\begin{proof}
Since $I$ is fully faithful, we have an isomorphism $\coll(I,I)\cong\coll(\id_{\A},\id_{\A})$ and so the inclusion  $J_0\colon\A\hookrightarrow\coll(\id_{\A},\id_{\A})$ of $\A$ in the collage as the first component has a right adjoint $R\colon\coll(\id_{\A},\id_{\A})\to\A$, that projects $\Ob(\A)\sqcup\Ob(\A)$ back to $\Ob(\A)$ and acts as the identity on homs. In particular, there is a $2$-cell $J_0R\Rightarrow\id_{\coll(\id_{\A},\id_{\A})}$ serving as the counit of the adjunction and the identity $2$-cell $\id_{\A}= RJ_0$ which is the unit. $\bbD$ sends this adjunction to a biadjunction $R^* \dashv_b J_0^*$ (they switch places because $\bbD$ is contravariant on $1$-cells) with a pseudonatural transformation $\alpha\colon R^*J_0^*\Rightarrow\id_{\bbD(\coll(\id_{\A},\id_{\A}))}$ as counit and a pseudonatural equivalence $\beta\colon\id_{\bbD(\A)}\simeq J_0^*R^*$ as unit.
Notice that in
\[
\begin{tikzcd}
J_1^*R^* \arrow[r, "J_1^*R^*\eta", Rightarrow] \arrow[d, "\simeq", sloped, phantom]  & J_1^*R^*I^*\Lan_I \arrow[r, "\simeq", Rightarrow] \arrow[d, "\simeq", sloped, phantom] & J_1^*R^*J_0^*\pi_{\B}^*\Lan_I \arrow[r, Rightarrow] \arrow[d, "\simeq", sloped, phantom] \arrow[rd, "\simeq", near start, xshift=-.5cm, yshift=.1cm, phantom] & J_1^*\pi_{\B}^*\Lan_I \\
\id_{\bbD(\A)} \arrow[r, "\eta"', Rightarrow] \arrow[ru, "\simeq", phantom] & I^*\Lan_I^* \arrow[r, "\simeq"', Rightarrow] \arrow[ru, "\simeq", phantom]      & J_0^*\pi_{\B}^*\Lan_I \arrow[ru, "\simeq"', sloped, Rightarrow]                                                        & {}                   
\end{tikzcd}
\]
the second arrows in the top and bottom rows are equivalences because $I^*\simeq J_0^*\pi_{\B}^*$, the vertical arrows are equivalences since $J_1^*R^*\simeq\id_{\bbD(\A)}$ and the oblique arrow at the right is an equivalence since $\pi_{\B}^*\simeq(IR)^*\simeq R^*I^*$ and finally $J_0^*\pi_{\B}^*\simeq J_0^*R^*I^*\simeq J_1^*R^*I^*\simeq J_1^*\pi_{\B}^*$.
Since the top row is an equivalence by the generalised elements version of \ref{axm: ptwise Kan ext}, we have that $\eta$ must be an equivalence as well.
\end{proof}
Similarly one can show that the counit $I^*\circ\Ran_I\Rightarrow\id_{\bbD(\A)}$ is a pseudonatural equivalence whenever $I$ is fully faithful.
\begin{prop}
If \ref{axm: ptwise Kan ext} holds for the $2$-prederivator $\bbD$, then the shifted 2-prederivator $\bbD^{\C}$ satisfies \ref{axm: ptwise Kan ext} for every $\C\in\Dia$.
\end{prop}
\begin{proof}
Let us prove the claim for left Kan extensions, as the proof for right Kan extensions is similar. We have to show that for every $f\colon\A\to\B$ in $\Dia$ and every $\C\in\Dia$, the diagram
\[\begin{tikzcd}
	\A & {\coll\B(f-,b)} \\
	\A & \B
	\arrow[from=1-1, to=1-2]
	\arrow[Rightarrow, no head, from=1-1, to=2-1]
	\arrow["{\pi_{\B}}", from=1-2, to=2-2]
	\arrow["f"', from=2-1, to=2-2]
\end{tikzcd}\]
is $\bbD^{\C}$-exact. If we consider the morphism $\C\times f\colon \C\times\A\to\C\times\B$, we know that   
\[\begin{tikzcd}
	\C\times\A & {\coll(\C\times f,\C\times b)} \\
	\C\times\A & \C\times\B
	\arrow[hook, from=1-1, to=1-2]
	\arrow[Rightarrow, no head, from=1-1, to=2-1]
	\arrow["{\pi_{\C\times\B}}", from=1-2, to=2-2]
	\arrow["{\C\times f}"', from=2-1, to=2-2]
\end{tikzcd}\]
is exact. The exactness of this diagram would imply the exactness of our original diagram, by the generalised version of \ref{axm: ptwise Kan ext}, as long as we show that $\coll(\C\times f,\C\times b)=\C\times\coll(f,b)$. This is indeed the case, given that they both have $(\C\times\A)\sqcup\C$ as set of objects and the morphisms between $\C\times\A$ and $\C$ are $(\C\times\B)((i, fa),(j, b))=\C(i,j)\times\B(fa,b)$. Therefore the two collages are the same and the claim is proven.
\end{proof}
\begin{prop}\label{HDer5 for rep pred}
If $\C$ is a complete and cocomplete $\sSet$-category, the represented $2$-prederivator $\Dia(-,\C)$ satisfies \ref{axm: ptwise Kan ext}.
\end{prop}
\begin{proof}
This comes from the expression of Kan extensions as weighted colimits in enriched category theory.
\end{proof}
\begin{prop}
\ref{axm: ptwise Kan ext} holds for the 2-prederivator $\bbD_{\M}$ associated to a combinatorial $\sSet_{\joy}$-enriched model category $\M$. 
\end{prop}
\begin{proof}
We need to prove that the mate of a certain 2-cell between homotopy categories is an equivalence.
We prove instead that the homotopy transformation induced by the mate of a point-set-level transformation is an equivalence.
The fact that these two coincide $-$ that is, that passing to homotopy commutes with the mate construction $-$ follows from the double pseudofunctor in Proposition \ref{hocat double pseudofunctor enriched}. Recall that Kan extensions of functors between categories enriched in simplicial sets can be computed using collages, therefore the square \[\begin{tikzcd}
	{[\A,\M]^\proj} && {[\B,\M]^\proj} && {[\B,\M]^\proj} & {[\B,\M]^\proj}  \\
	&&&& {[\coll(f,b),\M]^\proj} & {[\coll(f,b),\M]^\proj} \\
	{[\A,\M]^\proj} && {[\A,\M]^\proj} && {[\A,\M]^\proj} & {[\coll(f,b),\M]^\proj} 
	\arrow["{\Lan_f}", "\shortmid"{marking}, from=1-1, to=1-3]
	\arrow["\id", "\shortmid"{marking}, from=1-3, to=1-5]
	\arrow["\id", "\shortmid"{marking}, from=1-5, to=1-6]
	\arrow["\id"', from=1-1, to=3-1]
	\arrow["\id"', "\shortmid"{marking}, from=3-1, to=3-3]
	\arrow["{f^*}"', from=1-3, to=3-3]
	\arrow["{\pi_{\B}^*}"', from=1-5, to=2-5]
	\arrow["{i_{f,b}^*}"', from=2-5, to=3-5]
	\arrow["\id"', "\shortmid"{marking}, from=3-3, to=3-5]
	\arrow["\id"', "\shortmid"{marking}, from=2-5, to=2-6]
	\arrow["{\Lan_{i_{f,b}}}"', "\shortmid"{marking}, from=3-5, to=3-6]
	\arrow["{\pi_{\B}^*}", from=1-6, to=2-6]
	\arrow["\id", from=2-6, to=3-6]
	\arrow["\eta", shorten <=40pt, shorten >=40pt, Rightarrow, from=3-1, to=1-3]
	\arrow["\id"{description}, draw=none, from=3-3, to=1-5]
	\arrow["\id"{description}, draw=none, from=2-5, to=1-6]
	\arrow["\epsilon", shorten <=23pt, shorten >=23pt, Rightarrow, from=3-5, to=2-6]
\end{tikzcd}\]
in $\sSet_{\joy}$-\underline{Model} defined in \Cref{sec: double pseud}, is a simplicial natural transformation that is a componentwise weak equivalence. Calling $\alpha$ the composite square 
\[\begin{tikzcd}[ampersand replacement=\&]
	{[\A,\M]^{\proj}} \& {[\B,\M]^{\proj}} \\
	{[\A,\M]^{\proj}}\& {[\coll(f,b),\M]^{\proj}} 
	\arrow["\Lan_f", "\shortmid"{marking}, from=1-1, to=1-2]
	\arrow[Rightarrow, no head, xshift=-3mm, from=1-1, to=2-1]
	\arrow["{\Lan_{i_{f,b}}}"', "\shortmid"{marking}, from=2-1, to=2-2]
	\arrow["{\pi_{\B}^*}", xshift=-3mm, from=1-2, to=2-2]
	\arrow["\alpha", shorten <=20pt, shorten >=20pt, Rightarrow, from=2-1, to=1-2]
\end{tikzcd}\]
we have that the image of $\alpha$ through $\widetilde{h_*}$ specialises to a simplicial natural transformation whose components are the top rows in the diagram 
\[\begin{tikzcd}
	{(\Lan_{i_{f,b}}X_c)_f} & {(\pi_{\B}^*(\Lan_fX)_f)_c} \\
	{\Lan_{i_{f,b}}X_c} & {\pi_{\B}^*(\Lan_fX)_f} \\
	{\Lan_{i_{f,b}}X} & {\pi_{\B}^*(\Lan_fX)}
	\arrow["{\alpha_X}"', from=3-1, to=3-2]
	\arrow["\pi_{\B}^*j_{\Lan_fX}"', from=3-2, to=2-2]
	\arrow["\Lan_{i_{f,b}}r_{X}"', from=2-1, to=3-1]
	\arrow[from=2-1, to=2-2]
	\arrow["j_{\Lan_{i_{f,b}}X_c}", "\sim"', hook, from=2-1, to=1-1]
	\arrow["r_{\pi_{\B}^*(\Lan_fX)_f}", "\sim"', two heads, from=1-2, to=2-2]
	\arrow["\overline{\alpha}_X", dashed, from=2-1, to=1-2]
	\arrow["{\widetilde{\alpha}_X}", dashed, from=1-1, to=1-2]
\end{tikzcd}\]
with $X$ being a cofibrant-fibrant object of $[\A,\M]^{\proj}$, $j_Y\colon Y\xhookrightarrow{\sim} Y_f$ the trivial cofibration from an object $Y$ to its fibrant replacement, and $r_Y\colon Y_c\xrightarrowdbl{\sim} Y$ the trivial fibration to an object $Y$ from its cofibrant replacement $Y_c$. By Ken Brown's lemma $\Lan_{i_{f,b}}r_X$ is a weak equivalence because $r_X$ is a weak equivalence, both $X$ and $X_c$ are cofibrant objects, and $\Lan_{i_{f,b}}$ is left Quillen. Furthermore $\alpha_X$ is a weak equivalence by assumption and $\pi_{\B}^*j_{\Lan_fX}$ is a levelwise weak equivalence because $j_{\Lan_fX}$ is a projective trivial cofibration and therefore a levelwise trivial cofibration by Corollary \ref{equiv proj inj joyal enriched} so that $(\pi_{\B}^*j_{\Lan_fX})_a=(j_{\Lan_fX})_{\pi_{\B}(a)}$ is a trivial cofibration and in particular a weak equivalence for every $a\in\coll(f,b)$. Since weak equivalences in the projective model structure are levelwise, the composite $\pi_{\B}^*j_{\Lan_fX}\circ\alpha_X\circ\Lan_{i_{f,b}}r_X$ is a levelwise weak equivalence. Using the two-out-of-three property twice we conclude that $\widetilde{\alpha}_X$ is a levelwise weak equivalence between fibrant-cofibrant objects hence a levelwise homotopy equivalence\footnote{Notice that $\widetilde{\alpha}_X$ is the component at $X$ of a natural transformation between functors $[\A,\M]_{cf}^{\proj}\rightarrow[\coll(f,b),\M]_{cf}^{\proj}$ and so is a morphism in the codomain $[\coll(f,b),\M]_{cf}^{\proj}$. Namely $\widetilde{\alpha}_X$ is a natural transformation between functors $\coll(f,b)\to\M$.}. Using Lemma \ref{componentwise equiv} we deduce that its image through $h_*$ is a levelwise equivalence in $\Gray$, i.e.\ a pseudonatural equivalence. Therefore, by definition, its equivalence class in $h_2\Gray$ is an isomorphism which then lifts to pseudonatural equivalence in $\Gray$ giving the exactness of the diagrams in \ref{axm: ptwise Kan ext} for the prederivator $\bbD_{\M}$ associated to the combinatorial $\sSet_{\joy}$-model category $\M$.
\end{proof}

\subsubsection{Axiom on strongness} The next axiom also involves underlying diagram 2-functors. Roughly speaking it ensures that these 2-functors preserve certain powers in some weak sense. To express this more explicitly we will make use of smothering 2-functors (see also \cite[\S 3.6]{RVbook}).
\begin{defn}
A functor $F\colon\A\to\B$ is called \emph{smothering} if it is surjective on objects, full and conservative. 
\end{defn}
\begin{defn}
    A 2-functor is called \emph{smothering} if it is surjective on objects and it is locally smothering; in other words it is full on 1-cells and on 2-cells, and conservative on 2-cells. 
\end{defn}
\begin{mdframed}
    \begin{axiom}[(HDer 5)]\label{axm: strongness}
    The following $2$-functors are smothering for every $\A\in\Dia$: 
    \begin{enumerate}
        \item $\emph{dia}_{\Adj}^{\A}\colon\bbD(\A\times \Adj)\to[\Adj, \bbD(\A)]_p$,
        \item $\emph{dia}_{\b2}^{\A}\colon\bbD(\A\times\b2)\to[\b2, \bbD(\A)]_p$ and 
        \item $\emph{dia}_{\bbI}^{\A}\colon\bbD(\A\times\bbI)\to[\bbI,\bbD(\A)]_p$.
    \end{enumerate}
    where $\Adj$ is the free living adjunction of \cite{schanuel1986free}, and $\b2, \ \bbI$ are, respectively, the walking arrow and the walking isomorphism seen as locally discrete 2-categories.
\end{axiom}
\end{mdframed}
\begin{prop}
    The represented 2-prederivator satisfies \ref{axm: strongness}.
\end{prop}
\begin{proof}
    In this case these 2-functors are actually isomorphisms because we are working with representable 2-functors. 
\end{proof}
\begin{prop}
    If a 2-prederivator $\bbD$ satisfies \ref{axm: strongness}, then the shifted 2-prederivator $\bbD^{\B}$ satisfies \ref{axm: strongness} for every $\B\in\Dia$.
\end{prop}
\begin{proof}
    This is by design. Let us show it for 1.\ since the other points are completely analogous. For any $\B\in\Dia$, we want to show that the underlying diagram 2-functor $$\bbD^{\B}(\A\times \Adj)\to[\Adj, \bbD^{\B}(\A)]_p$$ is smothering. Notice that the LHS is isomorphic to $\bbD((\A\times\B)\times\Adj)$ by the universal property of the categorical product and the RHS is equal to $[\Adj,\bbD(\A\times\B)]_p$ by definition of the shifted 2-prederivator. Since point 1.\ of \ref{axm: strongness} is assumed to hold for $\bbD$, we know that $\dia_{\Adj}^{\A}$ is smothering for every simplicially enriched category $\A\in\Dia$. Hence we can conclude by choosing our simplicially enriched category to be $\A\times\B$.
\end{proof}
Using Remark \ref{shift of 2-pder of an enriched model category} we can simplify \ref{axm: strongness} for the 2-prederivator associated to a combinatorial $\sSet_{\joy}$-enriched model category, by requiring the underlying diagram 2-functors to be smothering in the non-shifted case. The proof that these 2-functors are surjective on objects is then analogous to the result on internalization of homotopy coherent adjunctions available in \cite{riehl2016homotopy}. The other aspects of smothering-ness remain conjectural; for instance, it is not immediately clear how to construct a morphism of homotopy coherent adjunctions starting from a morphism of adjunctions in the homotopy 2-category, nor how to lift an isomorphism between such morphisms to a fully coherent isomorphism. This would involve studying the functoriality of the constructions leading to the proof of the main result of \cite{riehl2016homotopy}.  
\begin{conj}
    The 2-prederivator $\bbD_{\M}$ associated to a combinatorial $\sSet_{\emph{Joyal}}$-enriched model category $\M$ satisfies \ref{axm: strongness}. 
\end{conj}


\section{Conclusions and future work}\label{chap:conclusions}
\graphicspath{{ch_other1/figures/}} 
In this section, we briefly explore potential directions for future research.

\subsection{Other proposed categories of diagrams}
In the current state of work, a 2-derivator is a weak kind of 3-functor. If we want to follow derivator theory more closely, we might need to consider strict 3-functors instead. A possible choice for the higher category of diagrams $\Dia$ might then be some 3-category of cofibrant simplicial categories (also known as simplicial computads) or even only cofibrant 2-categories, which are those whose underlying 1-category is the free category on a quiver. One could then be able to strictify 2-derivators by using similar constructions to the ones available in \cite{lack2002codescent, power1989general}.
\subsection{Axioms that still need to be checked}
\begin{enumerate}
\item In the current formulation of the theory of 2-prederivators, there is not a straightforward way to generalize the axiom (Der 2) so that the underlying diagram 2-functors become conservative for the represented 2-prederivator. In fact, the 2-cells in $\Dia$, which are the morphisms in the 2-category $\Dia(C,D)$ appearing when we are discussing a represented 2-prederivator \emph{are too strict}. A possible solution could be to replace the definition of a represented 2-prederivator with one based on the Gray hom $\GRAY(-,\K)$, and require that the axiom holds only on the subcategory of $\Dia$ spanned by 2-categories $-$ or, more generally, use the hom $\GRAY(h_*-, h_*\K)$. 
\item We still need to give a full proof of the smothering-ness of the underlying diagram 2-functors for the 2-derivator $\bbD_{\M}$ associated to a combinatorial $\sSet_{\joy}$-enriched model category $\M$ as in \ref{axm: strongness}.
\end{enumerate}
\subsection{Alternative to current axioms}
    It is known $-$ for example, from Proposition 1.26 in Groth's paper \cite{groth2013derivators} $-$ that the axiom (Der 4) in ordinary derivator theory is equivalent to the base change axiom for Grothendieck (op)fibrations. We propose a generalization of \ref{axm: ptwise Kan ext} based on this idea, as follows.
    \bigskip
    \begin{mdframed}
\begin{axiom}[(HDer 4')]
The following diagram in $\Dia$
\[\begin{tikzcd}
	\A & \B \\
	\C & \D
	\arrow["f", hook, from=1-1, to=1-2]
	\arrow[from=1-1, to=2-1]
	\arrow[from=1-2, to=2-2]
	\arrow[from=2-1, to=2-2]
	\arrow["\lrcorner"{anchor=center, pos=0.125, rotate=180}, draw=none, from=2-2, to=1-1]
\end{tikzcd}\]
    is exact whenever $f$ is a codiscrete cofibration (in the sense of \cite{street1980fibrations}). 
\end{axiom} 
\end{mdframed}
\subsection{Further axioms}
\begin{enumerate}  
\item So far, we did not explore the category theory of $\infty$-categories within a 2-derivator. We propose the following axiom as a starting point for this exploration, inspired by Observation 3.3.20 in \cite{riehl20152}. This axiom is intended to express the weak 2-universal property of the internal comma object, which is used in $\infty$-cosmology to study limits, adjunctions and other categorical constructions involving $\infty$-categories. One should be able to do the same in the base of every 2-derivator.

\bigskip 
\begin{mdframed}
 \begin{axiom}[(HDer 6)] 
Let $\boxslash$ and \reflectbox{$\Lcorner$} be the $2$-categories realizing respectively the shape of a comma square and of a pullback. Then the underlying diagram $2$-functor $\bbD(\boxslash)\to\bbD(\bb1)^{\boxslash}$ restricted to the image of $$\Ran_{i_{\text{\reflectbox{$\Lcorner$}}}}\colon\bbD(\reflectbox{$\Lcorner$})\to\bbD(\boxslash),$$ where $i_{\text{\reflectbox{$\Lcorner$}}}$ is the inclusion of \reflectbox{$\Lcorner$} in the square, induces a smothering functor $$\bbD(\bb1)(Y,f\downarrow g)\to\bbD(\bb1)(Y,f)\downarrow\bbD(\bb1)(Y,g),$$ where $Y\in\bbD(\bb1)$, $f$ and $g$ are the arrows in the diagram shape of the pullback, $f\downarrow g$ is the internal comma object and the target of the smothering functor is a comma object in $\twoCat$. 
\end{axiom}
\end{mdframed}

\bigskip
\item In the theory of derivators, a privileged place belongs to \emph{stable derivators}. These are an enhancement of triangulated categories that serve as an alternative to \emph{stable $\infty$-categories}. It would be nice to develop a similar approach to stable $(\infty,2)$-categories using stable 2-derivators, once a definition of stable 2-derivator becomes available. In order to do so, one should first identify the axioms that such a stable 2-derivator should satisfy.
\end{enumerate}
\subsection{Conjectures}
\begin{enumerate}
    \item The 2-prederivator $\bbD_{\M}$ can also be defined in terms of the injective model structure. We conjecture that this gives rise to a 2-prederivator equivalent to the one defined in this paper using the projective model structure, in the sense that there exists an equivalence of trihomomorphisms between them.
    \item In ordinary derivator theory, a result of Cisinski ensures that the homotopy theory of simplicial sets can be recovered in terms of derivators as the free completion of the point by homotopy colimits. One could try to do the same with 2-derivators to capture the category theory of $\infty$-categories.
    \item As mentioned in the \Cref{intro: domain} of the Introduction, one of the main motivations behind enlarging the domain in which 2-prederivators are defined is to capture monadicity in this context. In fact, from \cite{riehl2015completeness}, the forgetful functor from the category of algebras of an homotopy coherent monad can be built using weighted limits, as it is induced by a morphism of weights constructed via pushouts involving inclusions of the form $\partial\Delta^n\hookrightarrow\Delta^n$. In the language of 2-derivators, we can express weighted limits as follows: given a weight $W\colon\A\to\sSet$, the weighted limit 2-functor $\lim(W,-)$ is given by the composite
    \[\bbD(\A)\xrightarrow{\Ran_{i_{\A}^{W}}}\bbD(\coll W)\xrightarrow{\bullet^*}\bbD(\bb1),\]
    where $\bullet\colon\bb1\to\coll W$ picks the extra point in the collage. Thus, we want to be able to talk about diagram shapes that can be more general simplicial sets not arising as nerves of categories.
    \end{enumerate}
  
\section*{Acknowledgements}


This paper is based on work developed as part of my Ph.D.\ thesis \cite{NDVphd}. I would like to express my deepest gratitude to my main supervisor, Dominic Verity, for his invaluable help and for the many insightful discussions throughout the development of this work. I am also thankful to my associate supervisor, Steve Lack, for his careful proofreading and helpful suggestions.

I am grateful to the CoACT for providing a stimulating research environment during my doctoral studies. I would also like to thank Gabriele Campeggio and Pippo Baudo for their support. Additionally, I appreciate the anonymous referee for their valuable feedback, which helped improve the presentation of this paper.

My Ph.D.\ research was supported by an international Macquarie University Research Excellence Scholarship (iMQRES). I further acknowledge a Lift-Off Fellowship from the Australian Mathematical Society (AustMS), which provided financial support for my living expenses during the final stages of writing. Part of this article was presented during a visit to RIMS in Kyoto; I thank RIMS for their hospitality and the chance to discuss this topic.

Finally, I am deeply grateful to my parents for their unwavering moral support throughout my doctoral journey.

\printbibliography

@article{power1989general,
  title={A general coherence result},
  author={Power, A John},
  journal={Journal of Pure and Applied Algebra},
  volume={57},
  number={2},
  pages={165--173},
  year={1989},
  publisher={Elsevier}
}

@article{lack2002codescent,
  title={Codescent objects and coherence},
  author={Lack, Stephen},
  journal={Journal of Pure and Applied Algebra},
  volume={175},
  number={1-3},
  pages={223--241},
  year={2002},
  publisher={Elsevier}
}

@article{riehl2015completeness,
  title={Completeness results for quasi-categories of algebras, homotopy limits, and related general constructions},
  author={Riehl, Emily and Verity, Dominic},
  journal={Homology, Homotopy and Applications},
  volume={17},
  number={1},
  pages={1--33},
  year={2015},
  publisher={International Press of Boston}
}

@article{street1988gray,
  title={Gray’s tensor product of 2-categories},
  author={Street, Ross},
  journal={Manuscript, available at \url{https://web.science.mq.edu.au/~street/GrayTensor.pdf}},
  year={1988}
}

@InProceedings{Street:fib-yon,
author="Street, Ross",
editor="Kelly, Max",
title="Fibrations and Yoneda's lemma in a 2-category",
booktitle="Category Seminar",
year="1974",
OPTpublisher="Springer Berlin Heidelberg",
address="Berlin, Heidelberg",
pages="104--133",
isbn="978-3-540-37270-7"
}

@article{street1980fibrations,
  title={Fibrations in bicategories},
  author={Street, Ross},
  journal={Cahiers de {T}opologie et {G}{\'e}om{\'e}trie {D}iff{\'e}rentielle {C}at{\'e}goriques.},
  volume={21},
  number={2},
  pages={111--160},
  year={1980}
}

@article{StreetR:fortm,
title = "The formal theory of monads",
author = "Street, Ross",
year = "1972",
OPTmonth = "7",
volume = "2",
pages = "149--168",
journal = "Journal of Pure and Applied Algebra",
issn = "0022-4049",
publisher = "Elsevier",
number = "2",
}

@article{LackS:fortm,
title = "The formal theory of monads {II}",
author = "Lack, Steve and Street, Ross",
year = "2002",
OPTmonth = "11",
day = "8",
doi = "10.1016/S0022-4049(02)00137-8",
language = "English",
volume = "175",
pages = "243--265",
journal = "Journal of Pure and Applied Algebra",
issn = "0022-4049",
publisher = "Elsevier",
number = "1-3",
}

@phdthesis{cruttwell2008normed,
  title={Normed spaces and the change of base for enriched categories},
  author={Cruttwell, Geoff S. H.},
  year={2008},
  school={Dalhousie University}
}

@InProceedings{EilenbergS:cloc,
author = {Eilenberg, Samuel and Kelly, Max},
title = {Closed categories},
booktitle = {Proceedings of La Jolla Conference on Categorical Algebra},
pages = {421--562},
year = {1966},
publisher = {Springer-Verlag}
}

@Book{JohYau:2-cat,
author = {Johnson, Niles and Yau, Donald},
ALTeditor = {•},
title = {2-Dimensional Categories},
publisher = {Oxford University Press},
year = {2021},
OPTkey = {•},
OPTvolume = {•},
OPTnumber = {•},
OPTseries = {•},
OPTaddress = {•},
OPTedition = {•},
OPTmonth = {•},
OPTnote = {•},
OPTannote = {•}
}

@book{quillen2006homotopical,
  title={Homotopical algebra},
  author={Quillen, Daniel G.},
  volume={43},
  year={2006},
  publisher={Springer},
series = {Lecture Notes in Mathematics}
}

@book{hirschhorn2003model,
  title={Model categories and their localizations},
  author={Hirschhorn, Philip S.},
  number={99},
  year={2003},
  publisher={American Mathematical Society}
}

@article{schanuel1986free,
  title={The free adjunction},
  author={Schanuel, Stephen and Street, Ross},
  journal={Cahiers de {T}opologie et {G}{\'e}om{\'e}trie {D}iff{\'e}rentielle {C}at{\'e}goriques},
  volume={27},
  number={1},
  pages={81--83},
  year={1986}
}

@misc{NDVmres,
  author       = {Di Vittorio, Nicola}, 
  title        = {2-derivators},
  note         = {Master of Research Thesis, Macquarie University, \url{https://doi.org/10.25949/19817653.v1}, 2020}
  }

@misc{NDVphd,
  author       = {Di Vittorio, Nicola}, 
  title        = {Higher derivators as a foundation for $\infty$-category theory},
  note         = {Doctor of Philosophy Thesis, Macquarie University, \url{https://doi.org/10.25949/27885186.v1}, 2024}
  }

@book{grothendieck1983pursuing,
  title={Pursuing stacks},
  author={Grothendieck, Alexander},
  year={1983},
  publisher={unpublished manuscript}
}

@book{riehl2014categorical,
	title={Categorical {H}omotopy {T}heory},
	author={Riehl, Emily},
	number={24},
	year={2014},
	publisher={Cambridge University Press},
series = {New Mathematical Monographs}
}

@article{barwick2011unicity,
  author = {Barwick, Clark and Schommer-Pries, Christopher},
  title = {On the unicity of the homotopy theory of higher categories},
  journal = {arXiv},
  publisher = {Cornell University},
  url       = {https://arxiv.org/abs/1112.0040v4},
  pubstate      = {\bibstring{prepublished}},
  year = {2013}
}

@phdthesis{gurski2006algebraic,
  title={An algebraic theory of tricategories.},
  author={Gurski, Nick},
  year={2006},
  school={The University of Chicago}
}

@article{verity2011enriched,
  title={Enriched categories, internal categories and change of base},
  author={Verity, Dominic},
  journal={Reprints in {T}heory and {A}pplications of {C}ategories},
  number={20},
  pages={1--266},
  year={2011},
  publisher={Mount Allison University, Department of Mathematics and Science}
}

@book{kelly1982basic,
  title={Basic concepts of enriched category theory},
  author={Kelly, Max},
  volume={10},
  year={2005},
  publisher={Reprints in Theory and Applications of Categories}
}

@article{moser2019injective,
  title={Injective and projective model structures on enriched diagram categories},
  author={Moser, Lyne},
  journal={Homology, Homotopy and Applications},
  volume={21},
  number={2},
  pages={279--300},
  year={2019},
  publisher={International Press of Boston}
}

@article{shulman2011comparing,
  title={Comparing composites of left and right derived functors},
  author={Shulman, Michael},
  journal={New York J. Math},
  volume={17},
  pages={75--125},
  year={2011}
}

@article{bourke2017gray,
  title={The {G}ray tensor product via factorisation},
  author={Bourke, John and Gurski, Nick},
  journal={Applied Categorical Structures},
  volume={25},
  number={4},
  pages={603--624},
  year={2017},
  publisher={Springer}
}

@article{riehl20152,
  title={The 2-category theory of quasi-categories},
  author={Riehl, Emily and Verity, Dominic},
  journal={Advances in Mathematics},
  volume={280},
  pages={549--642},
  year={2015},
  publisher={Elsevier}
}

@article{groth2013derivators,
  title={Derivators, pointed derivators and stable derivators},
  author={Groth, Moritz},
  journal={Algebraic \& Geometric Topology},
  volume={13},
  number={1},
  pages={313--374},
  year={2013},
  publisher={Mathematical Sciences Publishers}
}

@misc{grothendieck1990derivateurs,
  title={Les d{\'e}rivateurs},
  author={Grothendieck, Alexander},
  year={1990}
}

@book{heller1988homotopy,
  title={Homotopy theories},
  author={Heller, Alex},
  volume={383},
  year={1988},
  publisher={American Mathematical Society}, 
  series = {Memoirs of the American Mathematical Society}
}

@book{lurie2009higher,
  title={Higher {T}opos {T}heory},
  author={Lurie, Jacob},
  year={2009},
  publisher={Princeton University Press}
}

@article{maltsiniotis2007k,
  title={La {K}-th{\'e}orie d'un d{\'e}rivateur triangul{\'e}.},
  author={Maltsiniotis, Georges},
  journal={Contemporary Mathematics},
  volume={431},
  pages={341--368},
  year={2007},
  publisher={Providence, RI: American Mathematical Society}
}

@article{maltsiniotis2001introduction,
  title={Introduction {\`a} la th{\'e}orie des d{\'e}rivateurs},
  author={Maltsiniotis, Georges},
  url={https://webusers.imj-prg.fr/~georges.maltsiniotis/ps/m.pdf},
  year={2001}
}

@article{more_elements,
  title={Elements of $\infty$-Category Theory},
  author={Riehl, Emily and Verity, Dominic},
  journal={\textit{Online draft with several additional chapters available at \url{https://elements-book.github.io/more-elements.pdf}}},
  year={2020}
}

@book{hovey2007model,
  title={Model categories},
  author={Hovey, Mark},
  number={63},
  year={2007},
  publisher={American Mathematical Society},
series = {Mathematical Surveys and Monographs}
}

@book{RVbook,
  title={Elements of $\infty$-Category Theory},
  author={Riehl, Emily and Verity, Dominic},
  volume={194},
  year={2022},
  publisher={Cambridge University Press},
series = {Cambridge Studies in Advanced Mathematics}
}

@book{boardman2006homotopy,
  title={Homotopy invariant algebraic structures on topological spaces},
  author={Boardman, John M. and Vogt, Rainer M.},
  volume={347},
  year={2006},
  publisher={Springer},
  series = {Lecture Notes in Mathematics}
}

@article{joyal2008theory,
  title={The theory of quasi-categories and its applications},
  author={Joyal, Andr{\'e}},  year={2008},
  journal={Quadern 45 vol II. Centre de Recerca Matem{\`a}tica Barcelona, \url{http://mat.uab.cat/~kock/crm/hocat/advanced-course/Quadern45-2.pdf},}
 }

@article{luriehigher,
  title={Higher {A}lgebra},
  author={Lurie, Jacob},
  year={2017},
  url={https://www.math.ias.edu/~lurie/papers/HA.pdf}
}

@phdthesis{lurie2004derived,
  title={Derived algebraic geometry},
  author={Lurie, Jacob},
  year={2004},
  school={Massachusetts Institute of Technology}
}

@article{gagna2022equivalence,
  title={On the equivalence of all models for $(\infty,2)$-categories},
  author={Gagna, Andrea and Harpaz, Yonatan and Lanari, Edoardo},
  journal={Journal of the London Mathematical Society},
  volume={106},
  number={3},
  pages={1920--1982},
  year={2022},
  publisher={Wiley Online Library}
}

@article{sava2022infty,
  author = {Sava, Chiara},
  title = {$\infty$-{D}old-{K}an correspondence via representation theory},
  journal = {arXiv},
  publisher = {Cornell University},
  url       = {https://arxiv.org/abs/2211.00762},
  pubstate      = {\bibstring{prepublished}},
  year = {2025}
}

@article{simpson1998homotopy,
  author = {Simpson, Carlos},
  title = {Homotopy types of strict 3-groupoids},
  journal = {arXiv},
  publisher = {Cornell University},
  url       = {https://arxiv.org/abs/math/9810059},
  pubstate      = {\bibstring{prepublished}},
  year = {1998}
}

@article{rezk2001model,
  title={A model for the homotopy theory of homotopy theory},
  author={Rezk, Charles},
  journal={Transactions of the American Mathematical Society},
  volume={353},
  number={3},
  pages={973--1007},
  year={2001}
}

@article{bergner2007model,
  title={A model category structure on the category of simplicial categories},
  author={Bergner, Julie},
  journal={Transactions of the American Mathematical Society},
  volume={359},
  number={5},
  pages={2043--2058},
  year={2007}
}

@book{gaitsgory2019study,
  title={A study in {D}erived {A}lgebraic {G}eometry: Volume I: {C}orrespondences and {D}uality},
  author={Gaitsgory, Dennis and Rozenblyum, Nick},
  volume={221.1},
  year={2017},
  publisher={American Mathematical Society},
  series = {Mathematical Surveys and Monographs}
}

@book{gaitsgory2017study,
  title={A study in {D}erived {A}lgebraic {G}eometry: Volume II: {D}eformations, {L}ie theory and formal geometry},
  author={Gaitsgory, Dennis and Rozenblyum, Nick},
  publisher={American Mathematical Society},
  volume={221.2},
series = {Mathematical Surveys and Monographs},
  year={2017}
}

@article{baez1995higher,
  title={Higher-dimensional algebra and topological quantum field theory},
  author={Baez, John C. and Dolan, James},
  journal={Journal of mathematical physics},
  volume={36},
  number={11},
  pages={6073--6105},
  year={1995},
  publisher={AIP Publishing}
}

@article{lurie2008classification,
  title={On the classification of topological field theories},
  author={Lurie, Jacob},
  journal={Current developments in mathematics},
  volume={2008},
  number={1},
  pages={129--280},
  year={2008},
  publisher={International Press of Boston}
}

@article{blumberg2013universal,
  title={A universal characterization of higher algebraic K-theory},
  author={Blumberg, Andrew J. and Gepner, David and Tabuada, Gon{\c{c}}alo},
  journal={Geometry \& Topology},
  volume={17},
  number={2},
  pages={733--838},
  year={2013},
  publisher={Mathematical Sciences Publishers}
}

@article{riehl2018homotopy,
author = {Riehl, Emily},
title = {Homotopy coherent structures},
journal = {Expositions in Theory and Applications of Categories,},
volume = {1},
number = {1},
pages = {1-31},
url = {http://www.tac.mta.ca/tac/expositions/articles/1/te1abs.html},
year = {2023}
}

@article{riehl2016homotopy,
  title={Homotopy coherent adjunctions and the formal theory of monads},
  author={Riehl, Emily and Verity, Dominic},
  journal={Advances in Mathematics},
  volume={286},
  pages={802--888},
  year={2016},
  publisher={Elsevier}
}

\end{document}